\documentclass[10pt]{article}
\usepackage[utf8]{inputenc}
\usepackage[english]{babel}
\usepackage{amsmath,amsthm,amsfonts}
\usepackage{graphicx}
\usepackage{bm}
\usepackage{tikz}
\usepackage{multirow}
\usepackage{color}

\graphicspath{ {./figs/} }

\title{A three-level multi-continua upscaling method for flow problems in fractured porous media}

\author{
Maria Vasilyeva \thanks{Institute for Scientific Computation, Texas A\&M University, College Station, TX 77843-3368 \& Department of Computational Technologies, North-Eastern Federal University, Yakutsk, Republic of Sakha (Yakutia), Russia, 677980. Email: {\tt vasilyevadotmdotv@gmail.com}.}
\and
Eric T. Chung \thanks{Department of Mathematics,
The Chinese University of Hong Kong (CUHK), Hong Kong SAR. Email: {\tt tschung@math.cuhk.edu.hk}.}
\and
Yalchin Efendiev \thanks{Department of Mathematics \& Institute for Scientific Computation (ISC),
Texas A\&M University,
College Station, Texas, USA. Email: {\tt efendiev@math.tamu.edu}.}
\and
Aleksey Tyrylgin 
\thanks{Multiscale model reduction laboratory, North-Eastern Federal University, Yakutsk, Republic of Sakha (Yakutia), Russia, 677980.}
}

\begin{document}

\maketitle

\begin{abstract}
Traditional two level upscaling techniques suffer from a high offline cost when the coarse grid size is much larger than the fine grid size. 
Thus, multilevel methods are desirable for problems with complex heterogeneities and high contrast. 
In this paper, we propose a novel three-level upscaling method for flow problems in fractured porous media. 
Our method starts with a fine grid discretization for the system involving fractured porous media. 
In the next step, based on the fine grid model,
we construct a nonlocal multi-continua upscaling (NLMC) method using an intermediate grid. The system resulting from NLMC
gives solutions that have physical meaning.
In order to enhance locality, the grid size of the intermediate grid needs to be relatively small, and this motivates using such an intermediate grid. 
However, the resulting NLMC upscaled system has a relatively large dimension. 
This motivates a further step of dimension reduction. In particular, we will apply the idea of the Generalized Multiscale Finite Element Method (GMsFEM)
to the NLMC system to obtain a final reduced model. 
We present simulation results for a two-dimensional model problem with a large number of fractures using the proposed three-level method.
\end{abstract}

\section{Introduction}

A fast and accurate solution of flow problems in fractured porous media is an important component in reservoir simulations. Direct numerical simulation requires using a very fine grid that resolves all scales and heterogeneities.
The resulting discrete formulation on the fine grid leads to a very large system of equations that is computationally expensive to solve.
To reduce the dimension of the system, multiscale methods or upscaling techniques are necessary \cite{houwu97, eh09, weinan2007heterogeneous, lunati2006multiscale, jenny2005adaptive}.
We will, in this paper, focus on a class of multiscale methods  based on local multiscale basis functions. 
In typical two level methods, multiscale basis functions are constructed locally, namely, within a coarse block or a union of several coarse blocks of an underlying coarse mesh,
which does not necessarily resolve any scale.
Constructing multiscale basis functions involves solutions, using the fine grid, of some local problems, which can be expensive for the case when coarse grid size is much larger than the fine grid size \cite{chung2016reiterated}. Therefore, problems with very large disparate scales require some coarsening techniques or multilevel techniques \cite{kunze2013multilevel}.
The commonly used techniques for such problems are the re-iterated homogenization methods or multilevel multiscale methods \cite{bensoussan2011asymptotic,lions2001reiterated, yuan2009hierarchical, szeliski1990fast, lipnikov2008multilevel, kunze2013multilevel, chung2016reiterated}. In multilevel multiscale approaches, multiple levels of coarsening are constructed by a recursive application of the basic two level method with the aim of improving computational efficiency. The main advantage of multilevel methods is to avoid solving local problems of large dimensions.

In our previous works, we developed multiscale model reduction techniques based on the Generalized Multiscale Finite Element Method (GMsFEM) for flow in fractured porous media \cite{akkutlu2015multiscale, chung2017coupling, efendiev2015hierarchical, akkutlu2018multiscale}. The general idea of GMsFEM is to design suitable spectral problems on some snapshot spaces to obtain dominant modes of the solutions. These dominant modes are used to construct the required multiscale basis functions \cite{EGG_MultiscaleMOR, egh12, chung2016adaptive, CELV2015}. The resulting multiscale space contains basis functions that take into account the microscale heterogeneities as well as high contrast and channelized effects, and the resulting multiscale scale solution  provides an accurate and efficient approximation of the fine scale solution.
We remark that the GMsFEM is related to the Proper Orthogonal Decomposition (POD) (c.f. \cite{egh12}) in the way that the GMsFEM constructs multiscale basis functions
that optimize an appropriate error within a finite dimensional space. The error of the GMsFEM has a spectral decay and is inversely proportional
to the eigenvalues of the spectral problems used for constructing basis functions.

Recently, the authors in \cite{chung2017constraint, chung2017non} proposed a new Constraint Energy Minimizing GMsFEM (CEM-GMsFEM)
with the aim of finding a multiscale method with a coarse mesh dependent convergence.
Constructing the multiscale space starts with an auxiliary space, which consists of eigenfunctions of a local spectral problem, and is defined for each coarse element. Using the auxiliary space, one can obtain the required multiscale basis functions by solving a constraint energy minimization problem. The resulting multiscale basis functions have an exponential decay away
from the coarse element for which the basis functions are formulated. Therefore, the multiscale basis functions are only numerically computed in an oversampled region defined by enlarging the target coarse element by a few coarse layers. It has been shown that these basis functions are able to capture high contrast channel effects. Moreover, the convergence of this method
depends only on the coarse grid size, and is independent of the scales and the heterogeneities of the coefficients of the PDE. We remark that the size of the oversampling domains depends on the coarse grid size and depends logarithmically on the contrast of the medium. 
Recently in \cite{chung2017non}, we introduced a non-local multi-continuum (NLMC) method for problems in heterogeneous fractured media. In the NLMC method, we construct multiscale basis functions based on the solution of some local constrained energy minimization problems as in the CEM-GMsFEM. 
One key ingredient of the NLMC method is that we can specify the location of all continua within coarse elements, and 
we construct these multiscale basis functions so that they have mean value zero in all continua within all coarse elements, except one target
continuum within a fixed coarse element. In this case, the degrees of freedoms of the resulting upscaled system have a physical meaning, namely,
they are the mean value of the solution on each continuum within each coarse element.
The NLMC has similar theoretical properties as that of the CEM-GMsFEM.


\begin{figure}[h!]
\centering
\includegraphics[width=0.6 \textwidth]{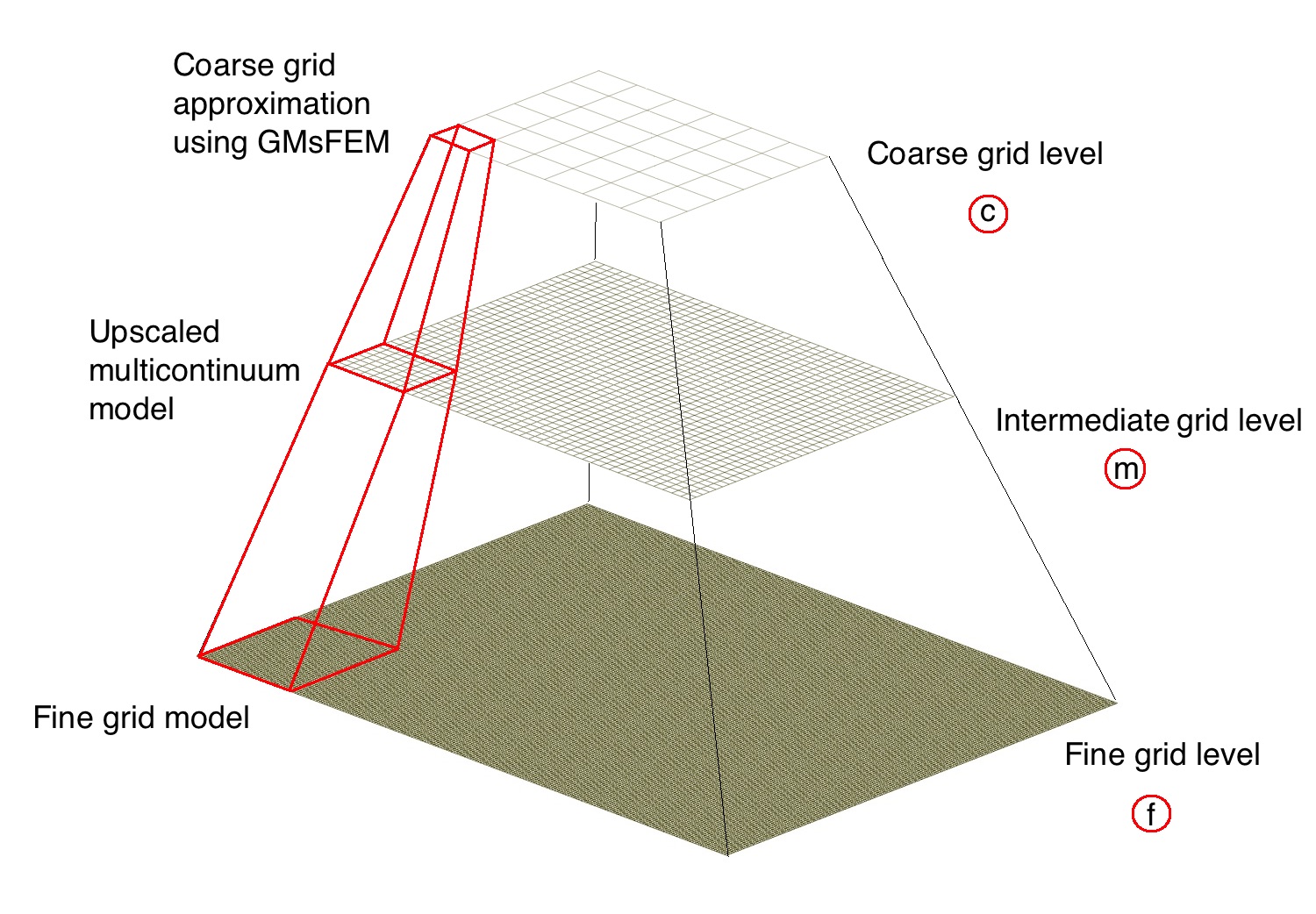}
\caption{Concept of three-level scheme.}
\label{fig:csch}
\end{figure}

As we mentioned above, two level multiscale methods can still suffer from large offline computational costs. 
In this work, we propose a new three level multiscale method based on both the GMsFEM and the NLMC
with the aim of taking advantage of both methodologies. Overall speaking, 
the proposed technique is the three-level scheme (see Figure \ref{fig:csch}) described as follows:
\begin{itemize}
\item fine grid model for fractured porous media,
\item intermediate grid model based on the NLMC method,
\item coarse grid approximation using the GMsFEM.
\end{itemize}
Our method starts with a fine grid discretization for the system involving fractured porous media. 
In the next step, based on the fine grid model,
we construct an NLMC method using an intermediate grid.  As discussed before, the system resulting from the NLMC method 
gives solutions that have physical meaning, namely, mean values on local continua. 
We remark that by an intermediate grid, we mean that the grid size is between the fine and the coarse grids. 
In order to enhance locality, the grid size of the intermediate grid needs to be relatively small, and this motivates
using such an intermediate grid. 
However, the resulting NLMC upscaled system has a relatively large dimension. 
This motivates a further step of dimension reduction. In particular, we will apply the idea of GMsFEM
to the NLMC system to obtain a final reduced model. 

This paper contains several novel ideas.
We present an extension of the GMsFEM for the NLMC models and show that the  GMsFEM can work with any multicontinuum upscaled model. The NLMC method provides an accurate upscaled multicontinuum approximation that we use for intermediate grid approximation.
The second advantage of the proposed method is the acceleration of the GMsFEM model construction, when  the solution of the local spectral problems are computationally expensive due to disparate scales and this requires coarsening \cite{chung2016reiterated, kunze2013multilevel}. Coarsening techniques should provide accurate and fast intermediate grid approximation. For this purpose, the NLMC method is applied for constructing the accurate upscaled intermediate grid model.

The paper is organized as follows. In Section 2, we consider a fine grid model to approximate the flow problem in the fractures porous media. In Section 3, we discuss an intermediate grid upscaled model construction using the NLMC method. Next, we present a construction of the multiscale basis functions on an intermediate grid for the GMsFEM in
Section 4 to obtain the final reduced model. Finally, we present numerical results and a conclusion in Section 5.

\section{Fine grid model}

First, we discuss the fine grid discretization of the flow system. 
We consider a mixed dimensional mathematical model for flow problem in fractured porous media.
A common approach to model fracture media is to consider the fractures as lower-dimensional objects \cite{martin2005modeling, d2012mixed, formaggia2014reduced, Quarteroni2008coupling}.
Let $\Omega \in \mathcal{R}^d$ (d = 2,3) be the computational domain for the porous medium and $\gamma \in \mathcal{R}^{d-1}$ be a reduced dimensional domain representing fracture networks. The flow model can be described as follows
\begin{equation}
\label{mm2}
\begin{split}
& a_m \frac{ \partial p_m }{\partial t} 
- \nabla \cdot (b_m \nabla p_m) +  \eta_m \sigma (p_m - p_f) =   q _m, 
\quad  x \in \Omega, \\
& a_f \frac{ \partial p_f }{\partial t} 
- \nabla \cdot (b_f \nabla p_f) -  \eta_f \sigma (p_m - p_f) =   q_f, 
\quad  x \in \gamma,
\end{split}
\end{equation}
\[
a_m = c_m, \quad a_f = d \, c_f, \quad 
b_m = k_m/\mu, \quad b_f = d \, k_f/\mu,
\]
where $\mu$  is the fluid viscosity, $c_{\alpha}$, $k_{\alpha}$  are the compressibility and permeability for porous matrix ($\alpha = m$) and fractured ($\alpha = f$), $q_{\alpha}$ is the source term for $\alpha = f, m$, $d$ is the fracture thickness, 
$p_m$ is the pressure in the porous matrix denoted by $\Omega$, 
$p_f$ is the pressure in the fractures $\gamma$. 
Coefficients $\eta_m$ and $\eta_f$ depend on mesh parameters and will be described later.

Let $\mathcal{T}_F=  \cup_i \varsigma_i$ be the fine grid with triangular or tetrahedral cells for the domain $\Omega$.
The  fracture mesh, denoted by $\mathcal{E}_{\gamma} = \cup_l \iota_l$, is constructed on the fractures domain $\gamma$. The coupled system (\ref{mm2}) is discretized using the embedded fracture model (EFM) \cite{hkj12, ctene2016algebraic, tene2016multiscale}. For the approximation in space, we apply the cell-centered finite-volume method with two-point flux approximation \cite{hkj12, ctene2016algebraic, bosma2017multiscale, ctene2017projection, tene2016multiscale}. 
Thus, we obtain the following discrete problem
\begin{equation}
\label{efm1}
\begin{split}
& a_m \frac{ p^{n+1}_{m, i} - p^{n}_{m, i} }{\tau} |\varsigma_i| 
 + \sum_j  T_{ij}  (p^{n+1}_{m, i} - p^{n+1}_{m, j})
 +  \sigma_{il} (p^{n+1}_{m, i} - p^{n+1}_{f, l} )  
 =  q_m   |\varsigma_i|, \quad \forall i = 1, N^m_F, \\
& a_f \frac{ p^{n+1}_{f, l} - p^{n}_{f, l}}{\tau}  |\iota_l| 
+ \sum_n W_{ln} (p^{n+1}_{f, l} - p^{n+1}_{f, n})
- \sigma_{il} (p^{n+1}_{m, i} - p^{n+1}_{f, l} ) 
 =  q_f  |\iota_l|, \quad \forall l = 1, N^f_F,
\end{split}
\end{equation}
where 
$T_{ij} = b_m |E_{ij}|/\Delta_{ij}$ ($|E_{ij}|$ is the length of facet between cells $\varsigma_i$ and $\varsigma_j$, $\Delta_{ij}$ is the distance between midpoint of cells $\varsigma_i$ and $\varsigma_j$),  
$W_{ln} = b_f/\Delta_{ln}$ ($\Delta_{ln}$ is the distance between points $l$ and $n$), 
$N^m_F$ is the number of cells in $\mathcal{T}_F$, 
$N^f_F$ is the number of cells related to the fracture mesh $\mathcal{E}_{\gamma}$, 
$\sigma_{il} = \sigma$ if $\iota_l \subset \varsigma_i$ and is zero otherwise. 
Here, we choose $\eta_m = 1 / |\varsigma_i|$, $\eta_f = 1 / |\iota_l|$ and use an implicit scheme for the time discretization, where $n$ is the number of time steps and $\tau$ is the given time step size. 

We can write the above scheme as the following system of equations for $p^{n} = (p^{n}_m, p^{n}_f)^T$ in matrix form
\begin{equation}
\label{mm-matrix}
M \frac{p^{n} -  p^{n-1} }{\tau} + A p^{n}= F,
\end{equation}
where
\[
M = 
\begin{pmatrix}
M_m & 0 \\
0 & M_f \\
\end{pmatrix}, \quad 
A = 
\begin{pmatrix}
A_m + Q & -Q \\
-Q & A_f+Q
\end{pmatrix}, \quad 
F  =
\begin{pmatrix}
F_m \\
F_f 
\end{pmatrix},
\]
and
\[
M_m = \{m^m_{ij}\}, \quad 
m^m_{ij} = 
\left\{\begin{matrix}
 a_m |\varsigma_i|  & i = j, \\ 
0 & i \neq j
\end{matrix}\right. , \quad 
M_f = \{m^f_{ln}\}, \quad 
m^f_{ln} = 
\left\{\begin{matrix}
 a_f |\iota_l|  & l = n, \\ 
0 & l \neq n
\end{matrix}\right. ,
\]\[
Q = \{q_{il}\}, \quad 
q_{il} = 
\left\{\begin{matrix}
\sigma & i = l, \\ 
0 & i \neq l
\end{matrix}\right. ,
\]
where
$A_m = \{T_{ij}\}$, 
$A_f = \{W_{ln}\}$, 
$F_m = \{f^m_i\}, \quad  f^m_i =  q_m |\varsigma_i|$, 
$F_f = \{f^f_l\}, \quad  f^m_i = q_f |\iota_l|$. We note that the size of this fine-grid system is $N_F = N^m_F + N^f_F$.

\section{The NLMC on intermediate grid}

In this section, we will construct an upscaled system for the fine system (\ref{mm-matrix}) on an intermediate grid. 
In particular, we will construct an upscaled model using the nonlocal multicontinua (NLMC) upscaling approach \cite{chung2017non}. 
In this method, the upscaled coefficients are based on the construction of multiscale basis functions. To do so, we solve local problems in some oversample local regions subject to the constraints that the mean values of the local solution vanishes in all continua except the one for which it is formulated. 
It has been shown that these multiscale basis functions have a spatial decay property and separate background medium and fractures. 
For more details in the derivation, we refer the reader to \cite{chung2017non}. Below, we will state a brief discussion of the derivation. 

Let $\mathcal{T}_I =  \cup_i K_i$ be a structured intermediate grid.
We consider a coarse cell $K_i$ and let $K_i^+$ be its oversampling region obtained by enlarging $K_i$ with few coarse cell layers.
For the fractures, we write $\gamma = \cup_{l = 1}^L \gamma^{(l)}$, where $\gamma^{(l)}$ denotes the $l$-th fracture network and $L$ is the total number of fracture networks.
Let $\gamma^{(l)}_j = K_j \cap \gamma^{(l)}$ be the fracture inside cell $K_j \in K_i^+$ and  $L_j$ be the number of fractures in $K_j$. 
For each $K_j \subset K_i^+$, we therefore need $L_j+1$ basis functions: one for $K_j$ and one for each $\gamma^{(l)}_j$. 
Following the framework of \cite{chung2017non} and \cite{chung2017constraint},
we will construct the required multiscale basis functions by solving a local problem on $K_i^+$
subject to some constraints to be specified in the following paragraph. 


We now define the constraints that will be used for multiscale basis construction. 
We use $\phi^{i,0}$ to denote the basis function corresponding to the porous matrix in the coarse element $K_i$
and use $\phi^{i,l}$ to denote the basis function corresponding to the $l$-th continuum within the coarse element $K_i$.
We remark that these basis functions are supported in $K_i^+$ and have zero trace on $\partial K_i^+$.
The required constraints are defined as follows: \\
(1) porous matrix in $K_i$, $\phi^{i,0} = (\phi^{i,0}_m, \phi^{i,0}_f)$ :
\[
\int_{K_j} \phi^{i,0}_m \, dx = \delta_{i,j}, \quad 
\int_{\gamma^{(l)}_j} \phi^{i,0}_f \, ds = 0, \quad  l=\overline{1, L_j}.
\]
(2) $l$-th fracture network in $K_i$, $\phi^{i,l} = (\phi^{i,l}_m, \phi^{i,l}_f)$:
\[
\int_{K_j} \phi^{i,l}_m \, dx = 0, \quad 
\int_{\gamma^{(l)}_j} \phi^{i,l}_f \, ds = \delta_{i,j}\delta_{m,l}, \quad  l=\overline{1, L_j}.
\]
We remark that the constraints are defined for each $K_j \subset K_i^+$.

To construct the multiscale basis functions with the energy minimizing property, we solve the following local problems in $K_i^+$ using a fine-grid approximation for flow in fractured porous media presented in  Section 2. In particular, we solve the following coupled system in  $K_i^+$:
\begin{equation}
\label{eq:basis}
\begin{pmatrix}
A^{i,+}_m+ Q^{i,+} & -Q^{i,+} & B^T_m & 0 \\
-Q^{i,+} & A^{i,+}_f + Q^{i,+} & 0 &  B^T_f \\
B_m & 0 & 0 & 0 \\
0 & B_f & 0 & 0 \\
\end{pmatrix} 
\begin{pmatrix}
\phi_m \\
\phi_f \\
\mu_m \\
\mu_f \\
\end{pmatrix} = 
\begin{pmatrix}
0 \\
0 \\
G_m \\
G_f \\
\end{pmatrix}
\end{equation}
with the zero Dirichlet boundary condition on $\partial K^+_i$ for both $\phi_m$ and $\phi_f$.  Here $A^{i,+}_m$, $A^{i,+}_f$ and $Q^{i,+}$ denote the parts of the fine-scale matrices that are related to the local domain $K^+_i$.  
Note that we used Lagrange multipliers $\mu_m$ and $\mu_f$ to impose the constraints defined above. 

For the construction of the multiscale basis function with respect to porous matrix $\phi^{i,0} = (\phi^{i,0}_m, \phi^{i,0}_f)$, we set $G_m = \delta_{i,j}$ and $G_f = 0$. For the multiscale basis function $\phi^{i,l} = (\phi^{i,l}_m, \phi^{i,l}_f)$ with respect to the $l$-th fracture network, we set  $G_m =  0$ and $G_f = \delta_{i,j}\delta_{m,l}$. Combining these multiscale basis functions, we obtain the following multiscale space
\[
V_{ms} = \text{span} \{ (\phi^{i,l}_m, \phi^{i,l}_f), \, i = \overline{1,N_c}, \, l = \overline{0, L_i} \}
\]
and the projection matrix
\[
R = \begin{pmatrix}
R_{mm} & R_{mf} \\
R_{fm} & R_{ff} \\
\end{pmatrix}, 
\]
where 
\[
R_{mm}^T = \left[ 
\phi^{0,0}_m, \phi^{1,0}_m  \ldots \phi^{N_c,0}_m
\right],\quad
 R_{ff}^T = \left[ 
\phi^{0,1}_f \ldots \phi^{0,L_0}_f, 
\phi^{1,1}_f \ldots \phi^{1,L_1}_f,
\ldots, 
\phi^{N_c,1}_f  \ldots \phi^{N_c,L_{N_c}}_f
 \right],
\]\[
R_{mf}^T = \left[ 
\phi^{0,0}_f, \phi^{1,0}_f \ldots \phi^{N_c,0}_f
 \right], \quad
R_{fm}^T = \left[ 
\phi^{0,1}_m \ldots \phi^{0,L_0}_m, 
\phi^{1,1}_m \ldots \phi^{1,L_1}_m,
\ldots, 
\phi^{N_c,1}_m  \ldots \phi^{N_c,L_{N_c}}_m
 \right], 
\]

Finally, the resulting upscaled intermediate grid model reads
\begin{equation}
\label{t-nlmc2}
\bar{M} \frac{\bar{p}^{n} - \bar{p}^{n-1}}{\tau} + \bar{A} \bar{p}^{n} = \bar{F},
\end{equation}
where 
$\bar{A} = R A R^T$, 
$\bar{p} = (\bar{p}_m, \bar{p}_f)$ is the average cell solution on intermediate grid element for porous matrix ($\bar{p}_m$) and for fractures ($\bar{p}_f$). We can reconstruct the downscale solution by $p = R^T \bar{p}$. 

As an approximation, we use diagonal mass matrix directly calculated on the intermediate grid 
\[
\bar{M} = 
\begin{pmatrix}
\bar{M}_m & 0 \\
0 & \bar{M}_f \\
\end{pmatrix}, \quad
\bar{F} = 
\begin{pmatrix}
\bar{F}_m \\
\bar{F}_f \\
\end{pmatrix}, 
\] 
where 
$\bar{M}_m = \text{diag}\{ a_m |K_i| \}$, 
$\bar{M}_f = \text{diag}\{ a_f |\gamma_i| \}$,  
and for the right-hand side vector 
$\bar{F}_m = \{ q_m |K_i| \}$, 
$\bar{F}_f = \{ q_f |\gamma_i| \}$. 
We remark that the matrix $A$ is non-local and provides a good approximation due to the coupling of various components in the basis construction.
The resulting upscaled model has one degree of freedom (DOF) for each fracture network and the size of intermediate grid system is $N_I = N^I_{cell} + \sum_{i=1}^{N^I_{cell}} L_i$, where $N^I_{cell}$ is the number of intermediate grid cells.

\section{The GMsFEM on coarse grid}

In this section, we will present a model reduction technique based on the GMsFEM. We will form a reduced model on a coarse grid 
based on the NLMC system constructed in the previous section. 
Generally speaking, the GMsFEM is a systematic approach to identify multiscale basis functions via local spectral problems \cite{egh12, EGG_MultiscaleMOR}. 
In the original GMsFEM, the method is constructed based on a fine grid discretization of the PDE.
In this paper, we will apply the GMsFEM idea to the system resulting from the NLMC method
and this is a new idea. 
To obtain a reduced system using GMsFEM, we first identify the local matrices from the NLMC system corresponding to
a set of overlapping coarse regions, typically called coarse neighborhoods \cite{egh12}. 
Then for each coarse neighborhood, we solve a spectral problem using the local matrices,
and select the dominant eigenfunctions corresponding to the small eigenvalues. 
The multiscale basis functions are then obtained by multiplying a suitable partition of unity function to the eigenfunctions. 
Finally, the GMsFEM system is obtained by forming a suitable projection matrix using the basis functions. 


For completeness, we summarize below the main steps in GMsFEM:
\begin{itemize}
\item[] \textit{Preprocessing (offline stage). } 
\begin{itemize}
\item The construction of the multiscale basis functions in local domains.
\item The construction of the coarse grid system.
\end{itemize}
\item[] \textit{Solver (online stage).} 
\begin{itemize}
\item Solution of the coarse grid system.
\end{itemize}
\item[] \textit{Postprocessing.} 
\begin{itemize}
\item Reconstruction of the fine grid solution.
\end{itemize}
\end{itemize}

In the following, 
we will describe the construction of the multiscale basis functions $\psi_k^{\omega}$ which is supported in a coarse neighborhood
$\omega$, where $k$ represents the numbering of the basis functions.  

Let $\mathcal{T}_C = \cup_i \Theta_i$ be the structured coarse grid and assume that each coarse element is a connected union of fine grid and intermediate grid blocks. 
We use $\{x_i\}_{i=1}^{N^C_{vert}}$ to denote the vertices of the coarse mesh $\mathcal{T}_C$, where $N^C_{vert}$ is the number of coarse nodes.  We define the coarse neighborhood of the node $x_i$ by $\omega_i = \cup_{j}\left\{ \Theta_j | \,  x_i \in \overline{\Theta}_j \right\}$. 

We now consider a coarse neighborhood $\omega_i$.
In order to construct the multiscale space $V_{ms}^{\omega_i}$ with respect to $\omega_i$, 
we solve following local spectral problem in local domain $\omega_i$
\begin{equation} 
\label{offeig}
A \Psi^i = \lambda^i S \Psi^i,
\end{equation}
where the matrix $A$ is the restriction of the matrix $\bar{A}$ in the coarse neighborhood $\omega_i$
and the matrix $\bar{A}$ is the matrix resulting from the NLMC method \eqref{t-nlmc2}. 
Moreover, the matrix $S$ is defined as follows:
\[
S =
\begin{pmatrix}
\bar{S}_m & 0 \\
0 & \bar{S}_f \\
\end{pmatrix},  \quad
\bar{S}_m = \{s^m_{ij}\}, \quad 
s^m_{ij} = 
\left\{\begin{matrix}
 b_m |K_i|  & i = j, \\ 
0 & i \neq j
\end{matrix}\right. , \quad 
\bar{S}_f = \{s^f_{ln}\}, \quad 
s^f_{ln} = 
\left\{\begin{matrix}
 b_f |\gamma_l|  & l = n, \\ 
0 & l \neq n
\end{matrix}\right. .
\] 
To define the required multiscale space, we choose eigenvectors $\Psi^i_k$ ($k=1,...,M_i$) corresponding to the smallest $M_i$ eigenvalues and set
\begin{equation} 
\label{cgspace}
V_C  = \text{span} \{ \psi^i_k = \chi_i  \Psi_k^i : \, 1 \leq i \leq N^C_{vert} \quad  \text{and} \quad  1 \leq k \leq M_i  \},
\end{equation}
where $\chi_i$ are the standard linear partition of unity functions and $M_i$ denotes the number of eigenvectors that are chosen for each coarse node $i$. 
The construction in \eqref{cgspace} yields a counterpart of the continuous basis functions due to the multiplication of local domain eigenvectors with the continuous partition of unity functions. 

Using a single index notation for the basis functions, we may write 
\[
V_C = \text{span} \{ \psi_1, \psi_2, ... , \psi_{N_C} \},  \quad
R_C^T = \left[ \psi_1 , \ldots, \psi_{N_C} \right],
\] 
where $R_C$ is the projection matrix and  $N_C =\sum_{i=1}^{N^C_{vert}} M_i$ is the  size of the coarse grid system.
Finally, we can write the GMsFEM system as 
\begin{equation}
\label{t-nlmc3}
M_C \frac{p^{n}_C  - p^{n-1}_C }{\tau} + A_C p^{n}_C = F_C,
\end{equation}
and $p_C \in V_C$ and $p_C = \sum_i p_{C,i} \psi_i(x)$.
In the above system, 
we have 
\[
M_C = R_C \bar{M} R_C^T, \quad 
A_C = R_C \bar{A} R_C^T, \quad 
F_C = R_C \bar{F},
\]
and $\bar{p} = R_C^T p_C$ is the reconstructed intermediate grid solution and $p = R^T \bar{p} $  is the reconstructed fine grid solution.

\section{Numerical results}

\begin{figure}[h!]
\centering
\includegraphics[width=0.49 \textwidth]{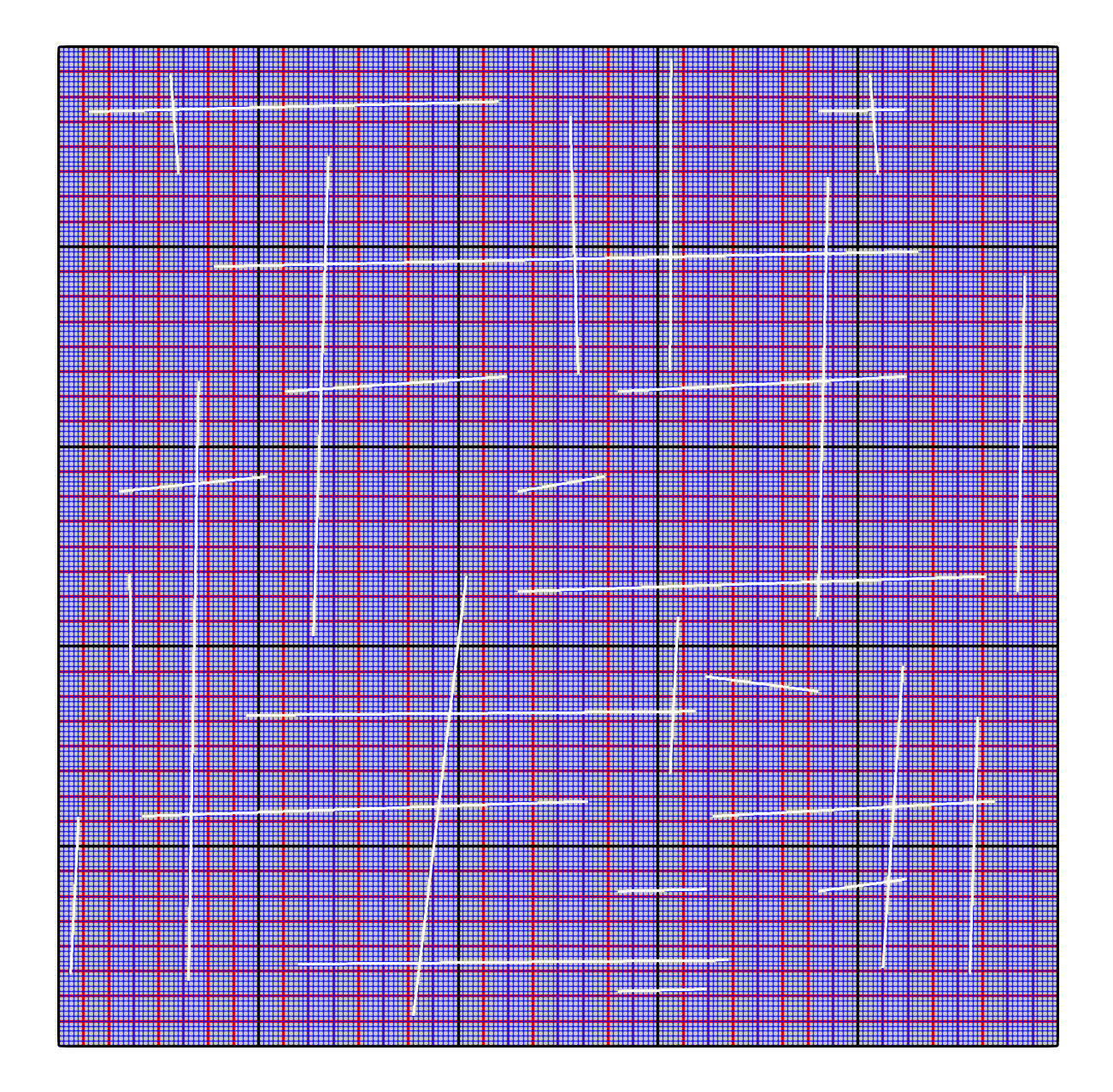}
\includegraphics[width=0.49 \textwidth]{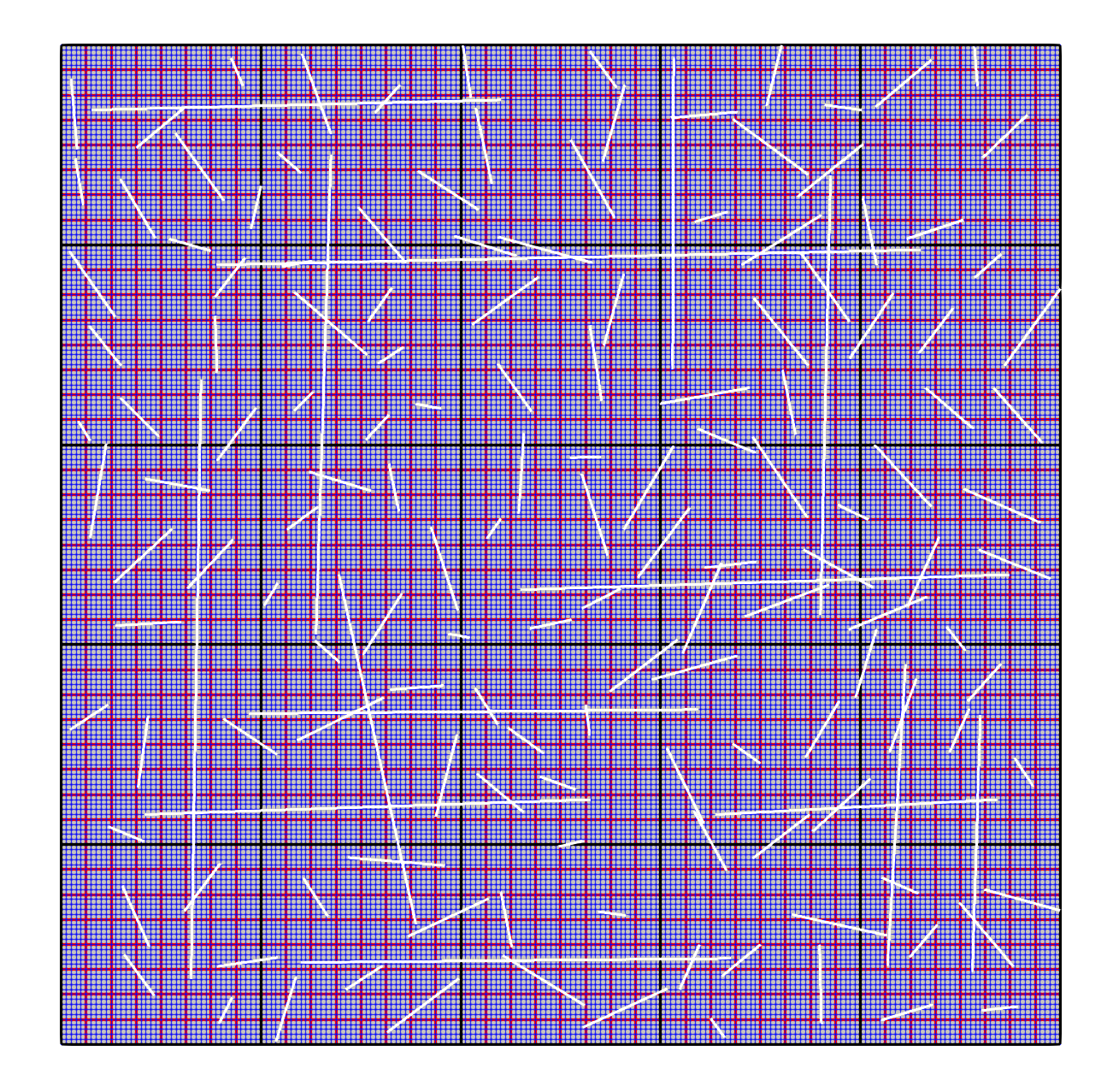}
\caption{Computational grids (black color - coarse grid, red color - intermediate grid and blue color - fine grid). Fractures are depicted by white color. 
Left: \textit{Geometry 1} with 30 fracture lines. 
Right: \textit{Geometry 2}  with 160 fracture lines. }
\label{fig:mesh}
\end{figure}

In this section, we present numerical results for our three level scheme. We consider the problem in domain $\Omega = [0, 1] \times [0, 1]$.  As model problems, we consider two geometries with different fracture distribution:
\begin{itemize}
\item \textit{Geometry 1}. Domain with 30 fracture lines. 
\item \textit{Geometry 2}. Domain with 160 fracture lines. 
\end{itemize}
In Figure \ref{fig:mesh}, we show computational grids for \textit{Geometry 1} and \textit{Geometry 2}. 
The implementation is based on the open-source simulation library FEniCS, where we use geometry objects and interface to the linear and spectral solvers \cite{logg2009efficient, logg2012automated}.

\begin{figure}[h!]
\centering
\includegraphics[width=0.32 \textwidth]{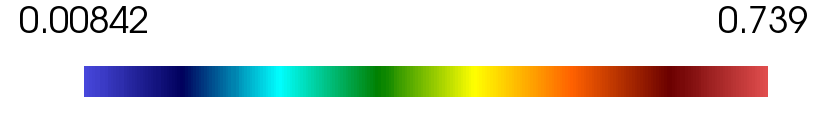}\\
\includegraphics[width=0.32 \textwidth]{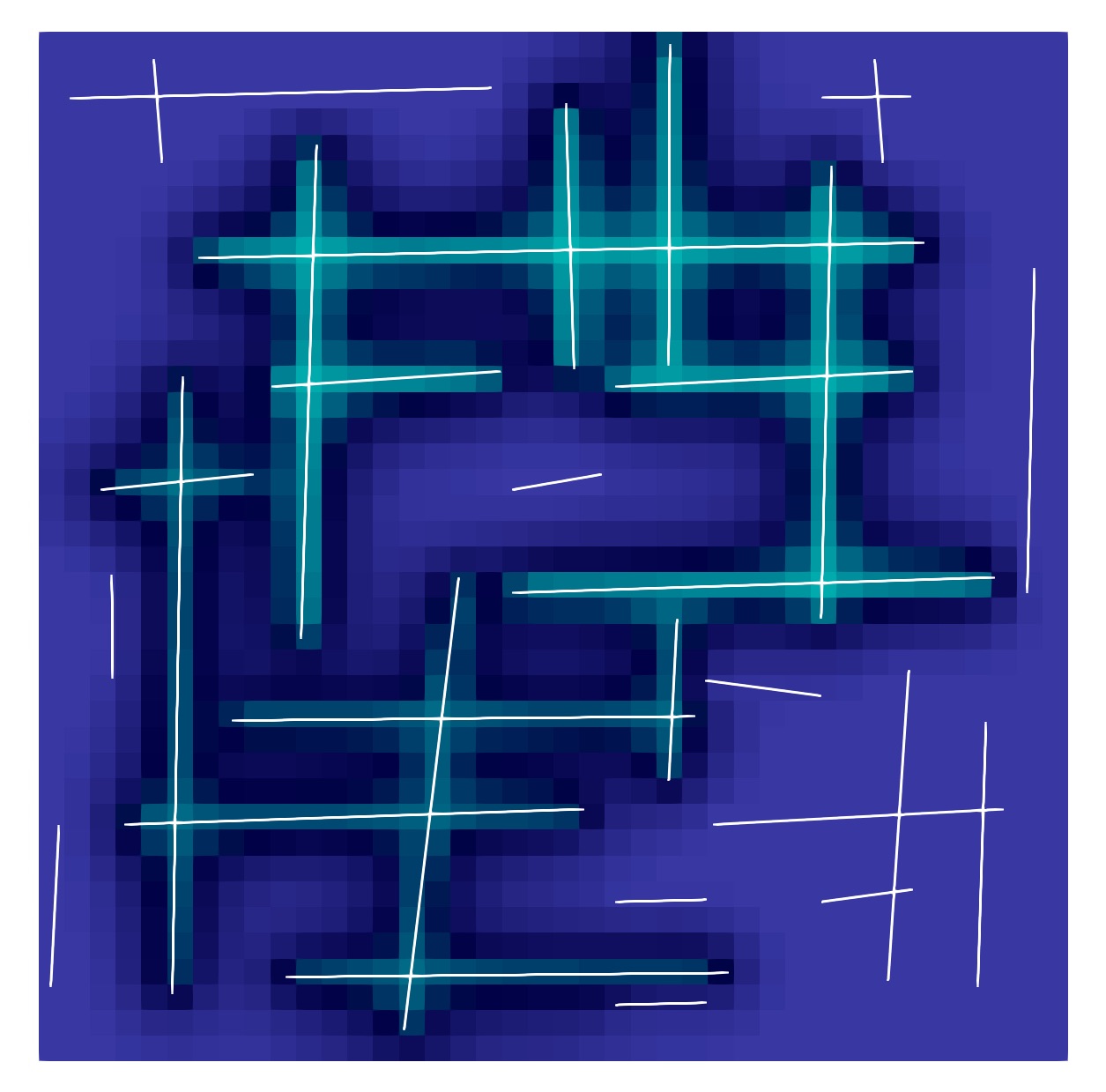}
\includegraphics[width=0.32 \textwidth]{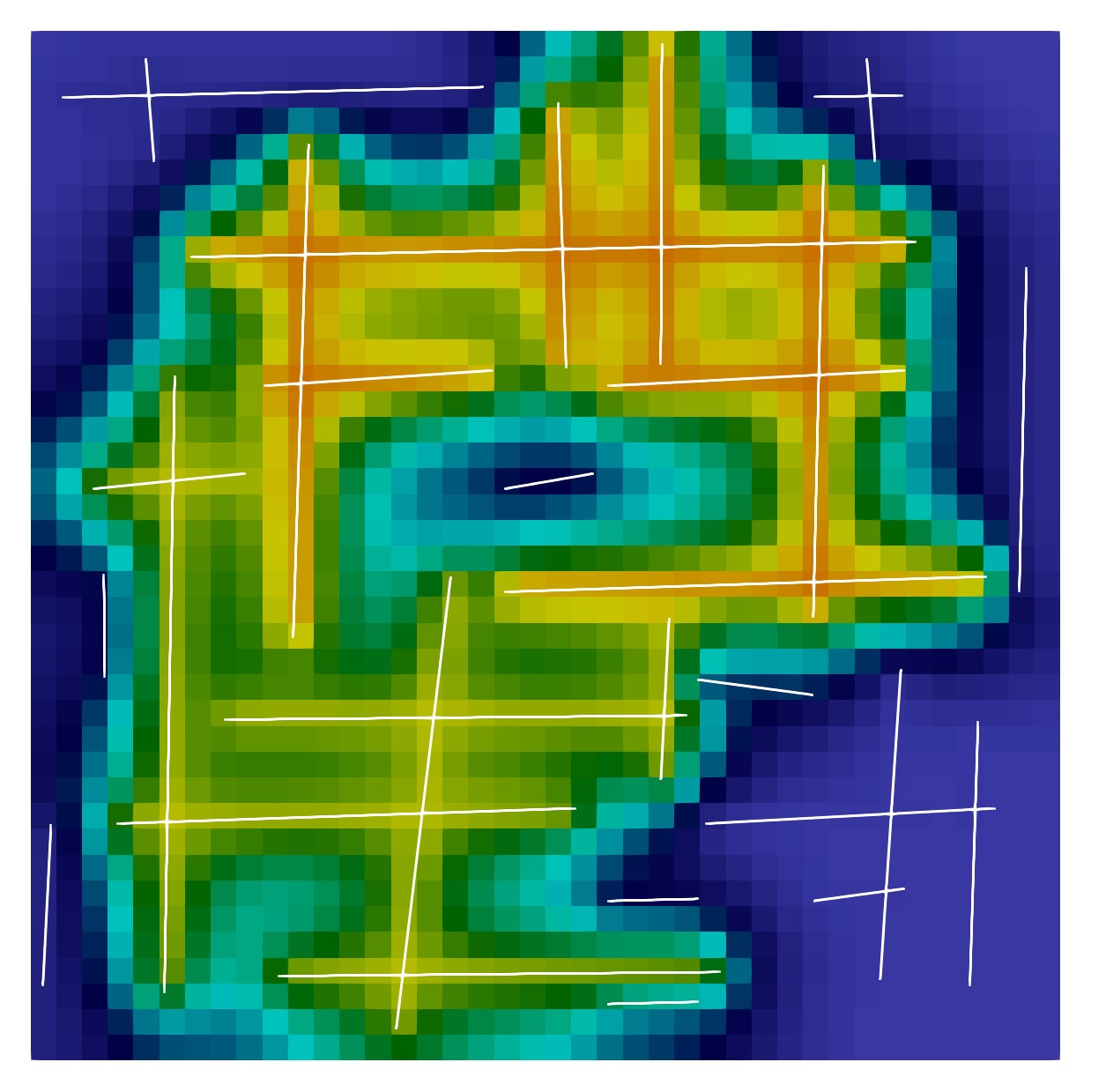}
\includegraphics[width=0.32 \textwidth]{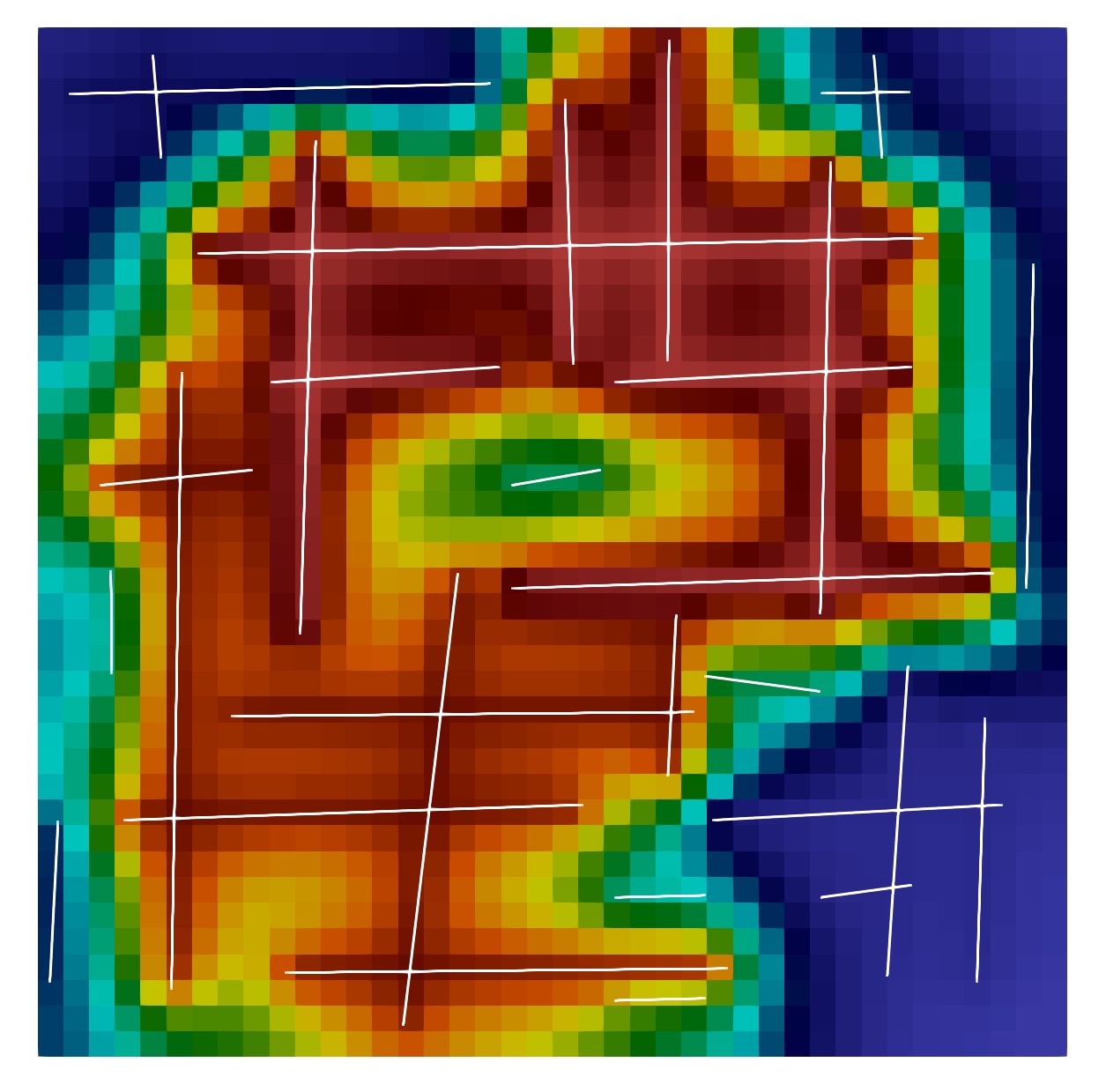}\\
\includegraphics[width=0.32 \textwidth]{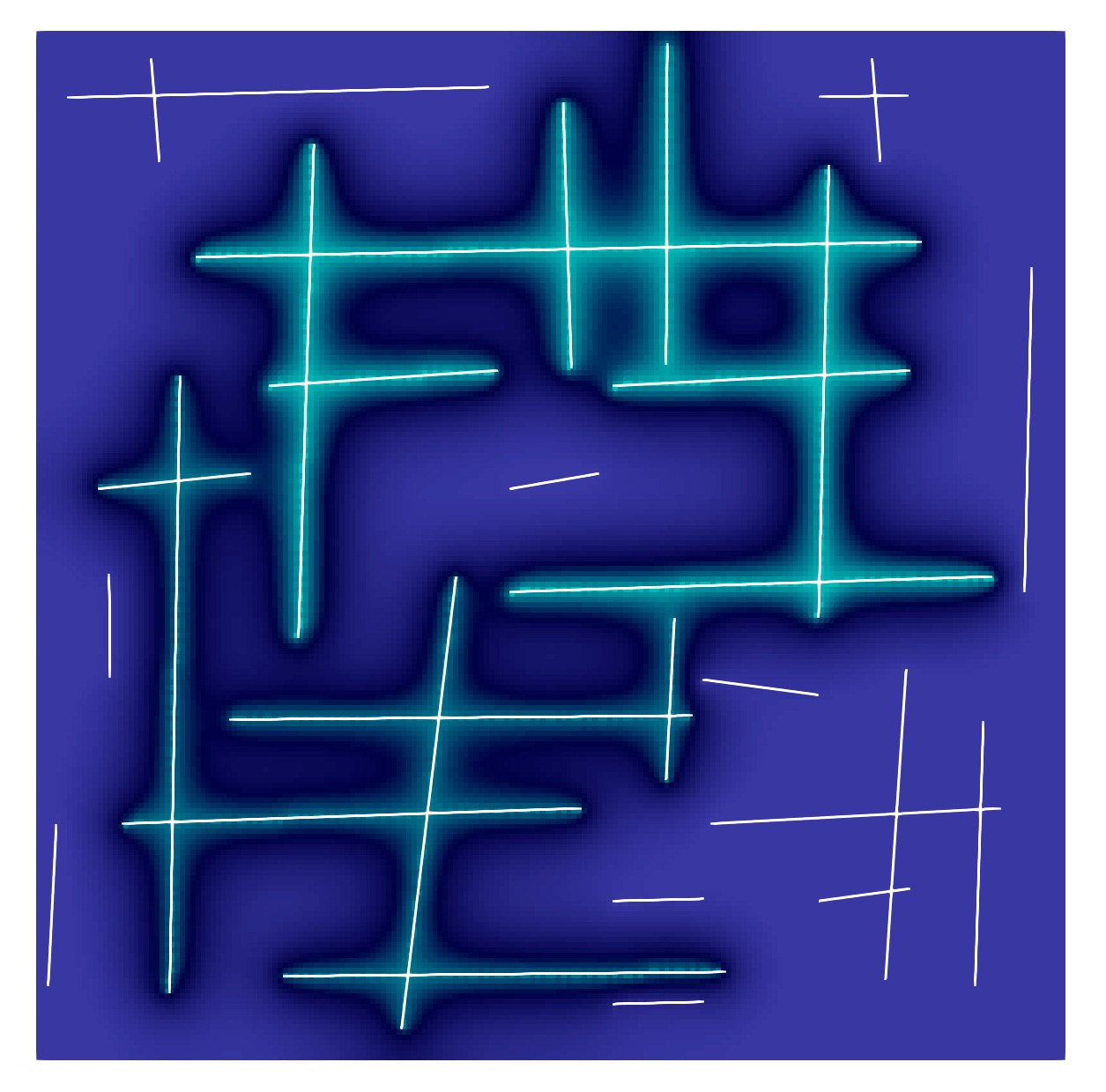}
\includegraphics[width=0.32 \textwidth]{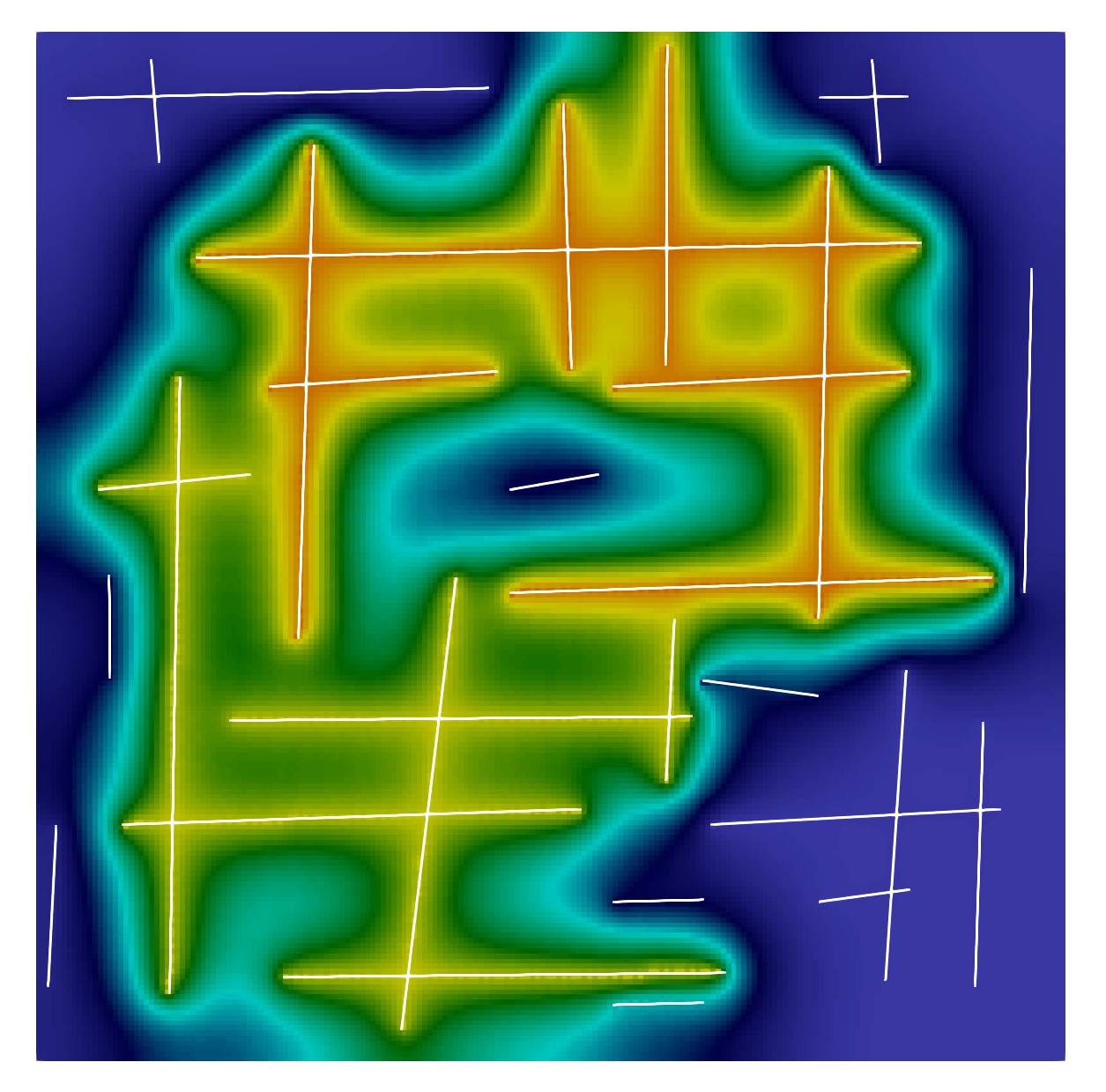}
\includegraphics[width=0.32 \textwidth]{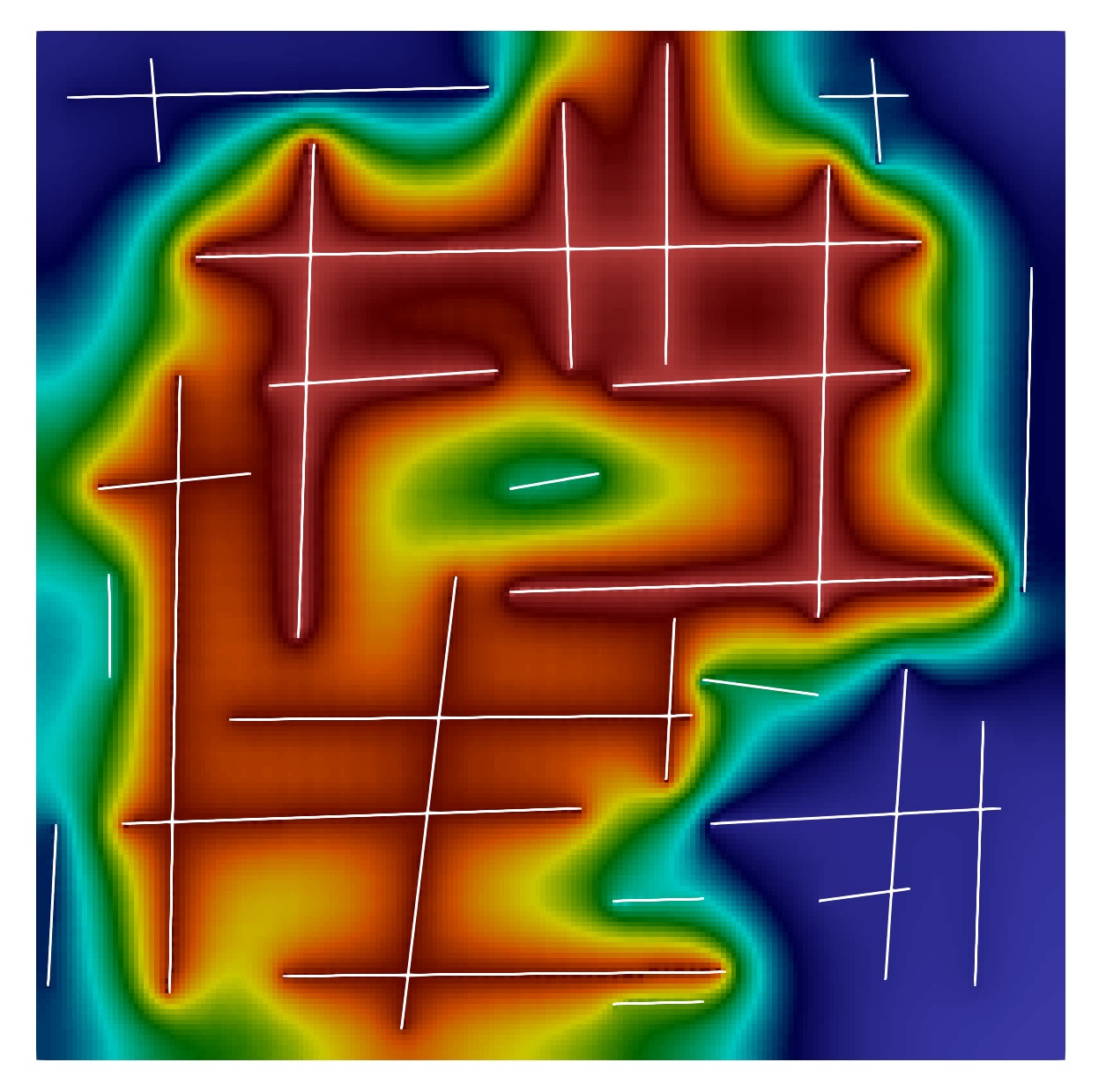}
\caption{Multiscale solutions on mesh $40 \times 40$ with $K^4$ using NLMC model for different time steps $t_{10} = 0.02$,  $t_{30} = 0.06$ and $t_{50} = 0.1$ (from top to bottom). \textit{Geometry 1}. 
First row: upscaled intermediate grid solution. 
Second row: downscaled fine grid solution.  }
\label{fig:u1}
\end{figure}

\begin{figure}[h!]
\centering
\includegraphics[width=0.32 \textwidth]{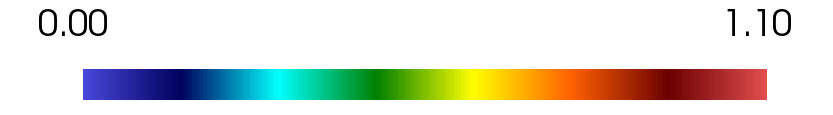}\\
\includegraphics[width=0.32 \textwidth]{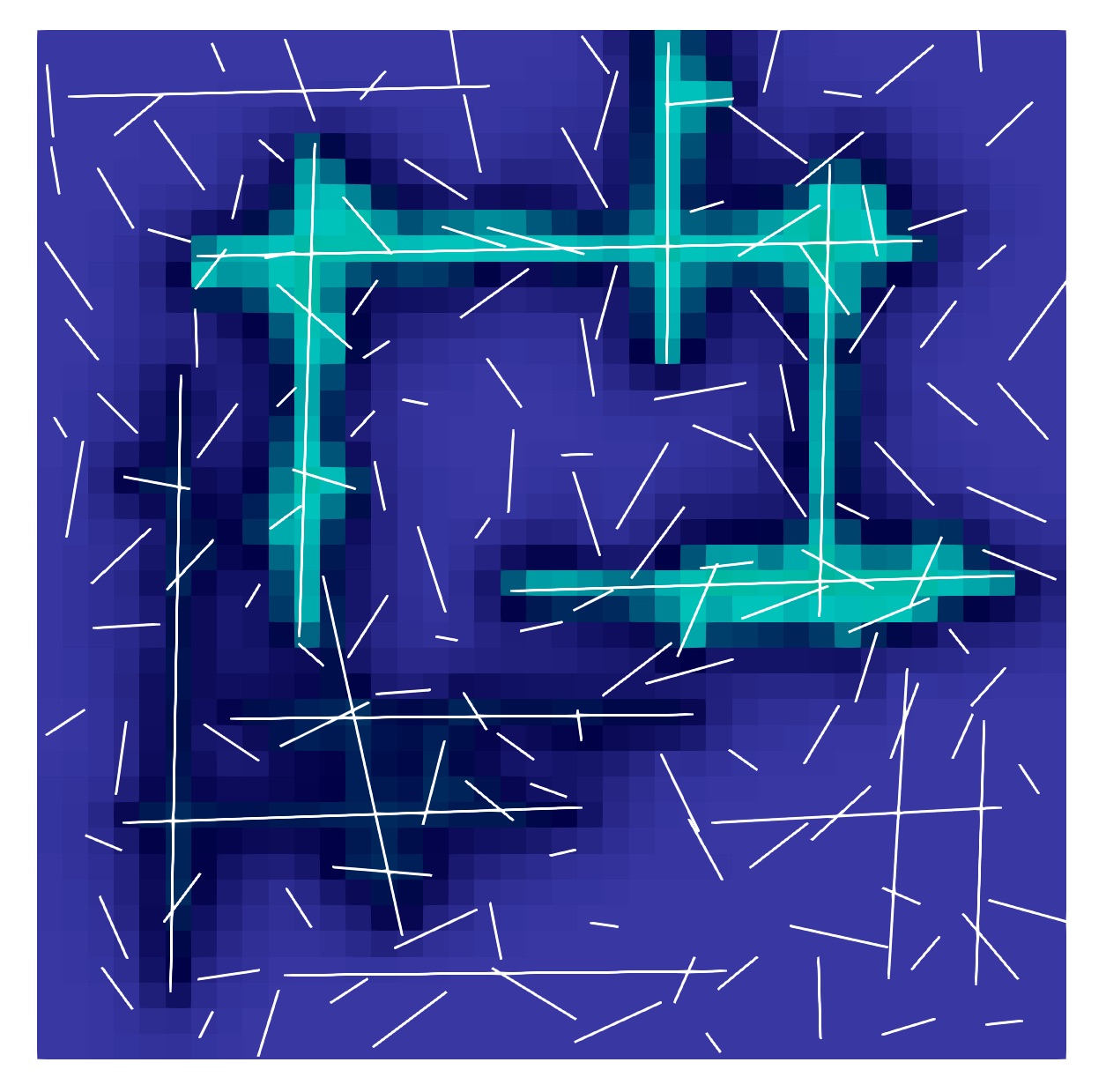}
\includegraphics[width=0.32 \textwidth]{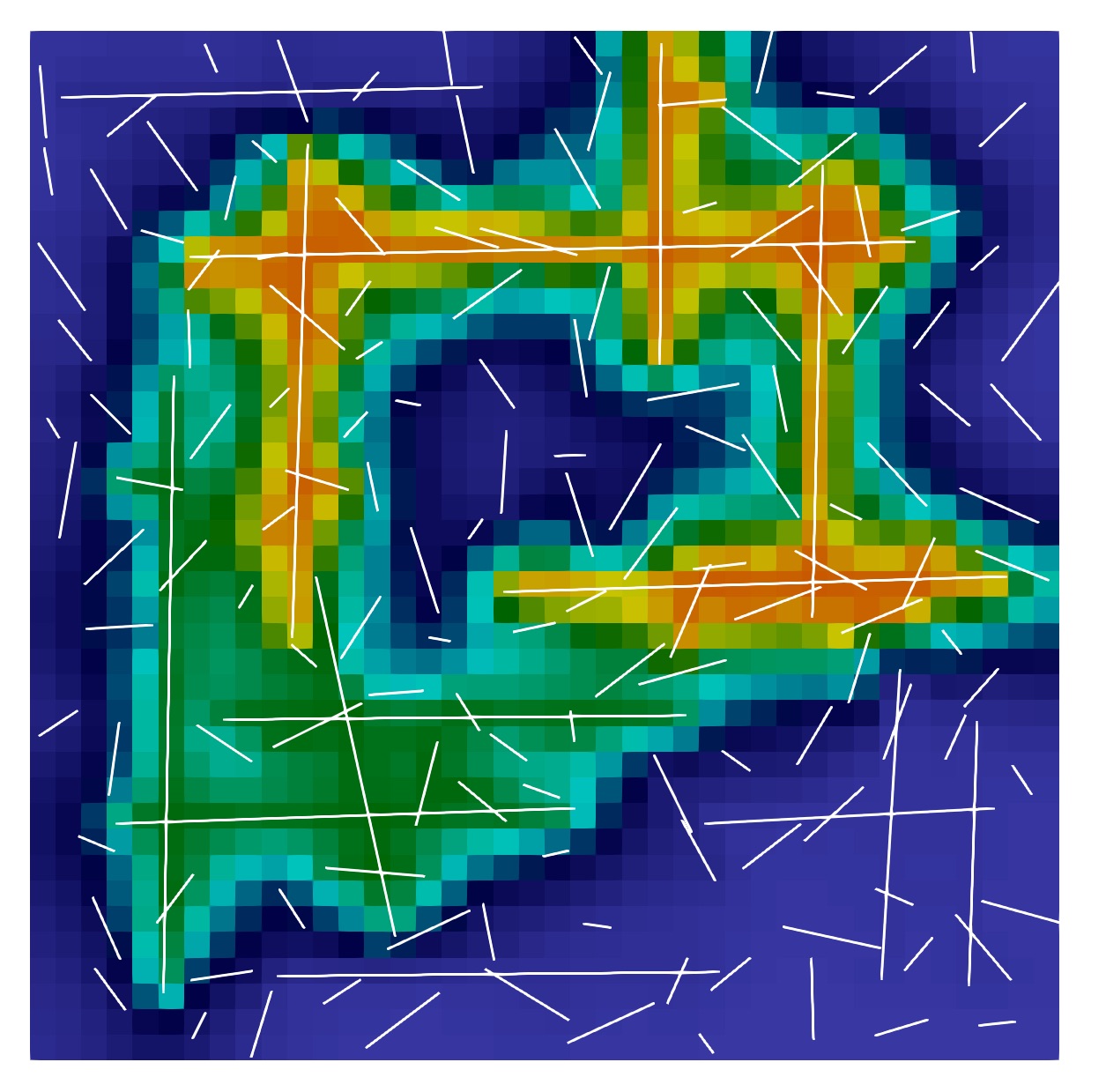}
\includegraphics[width=0.32 \textwidth]{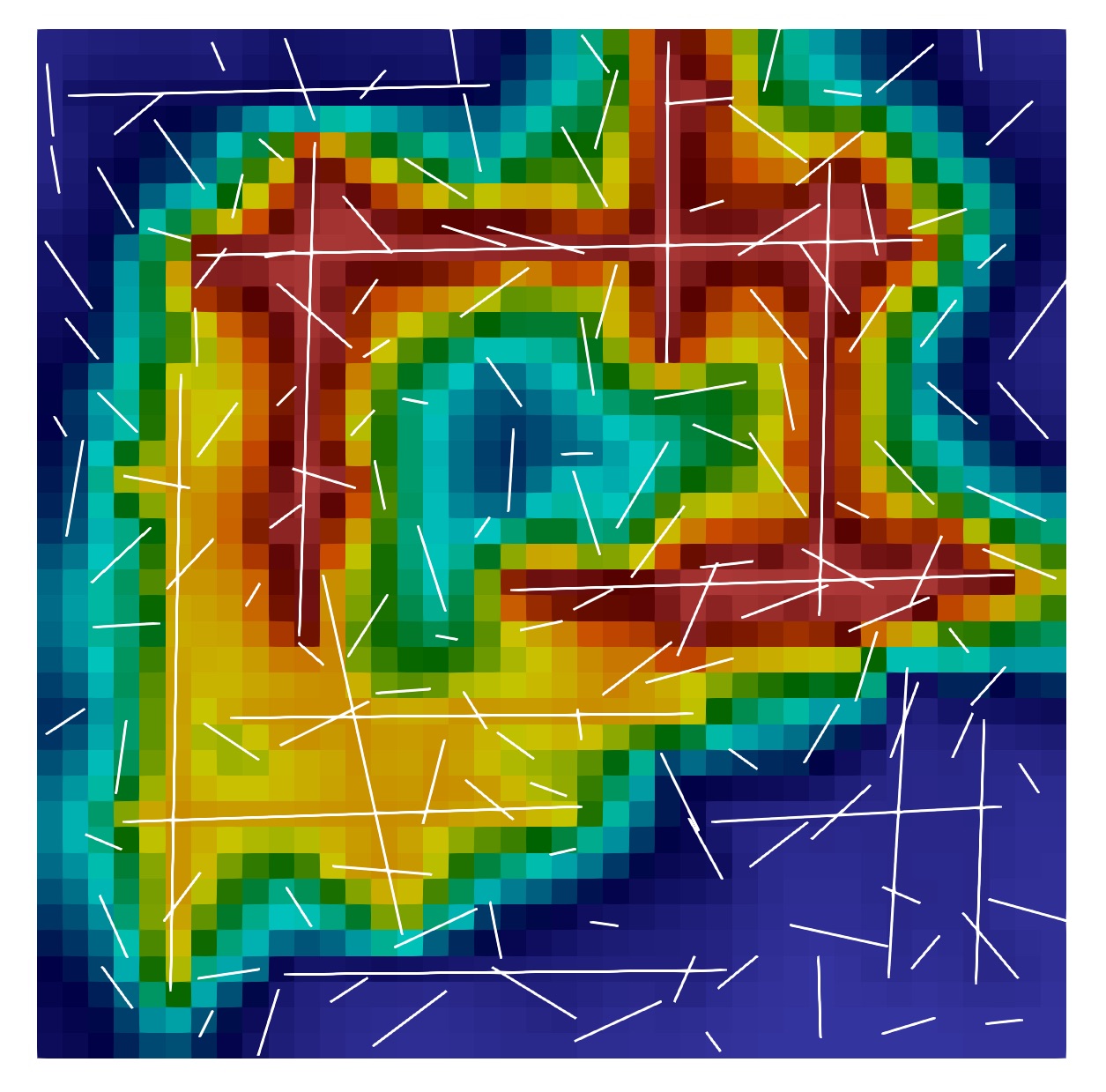}\\
\includegraphics[width=0.32 \textwidth]{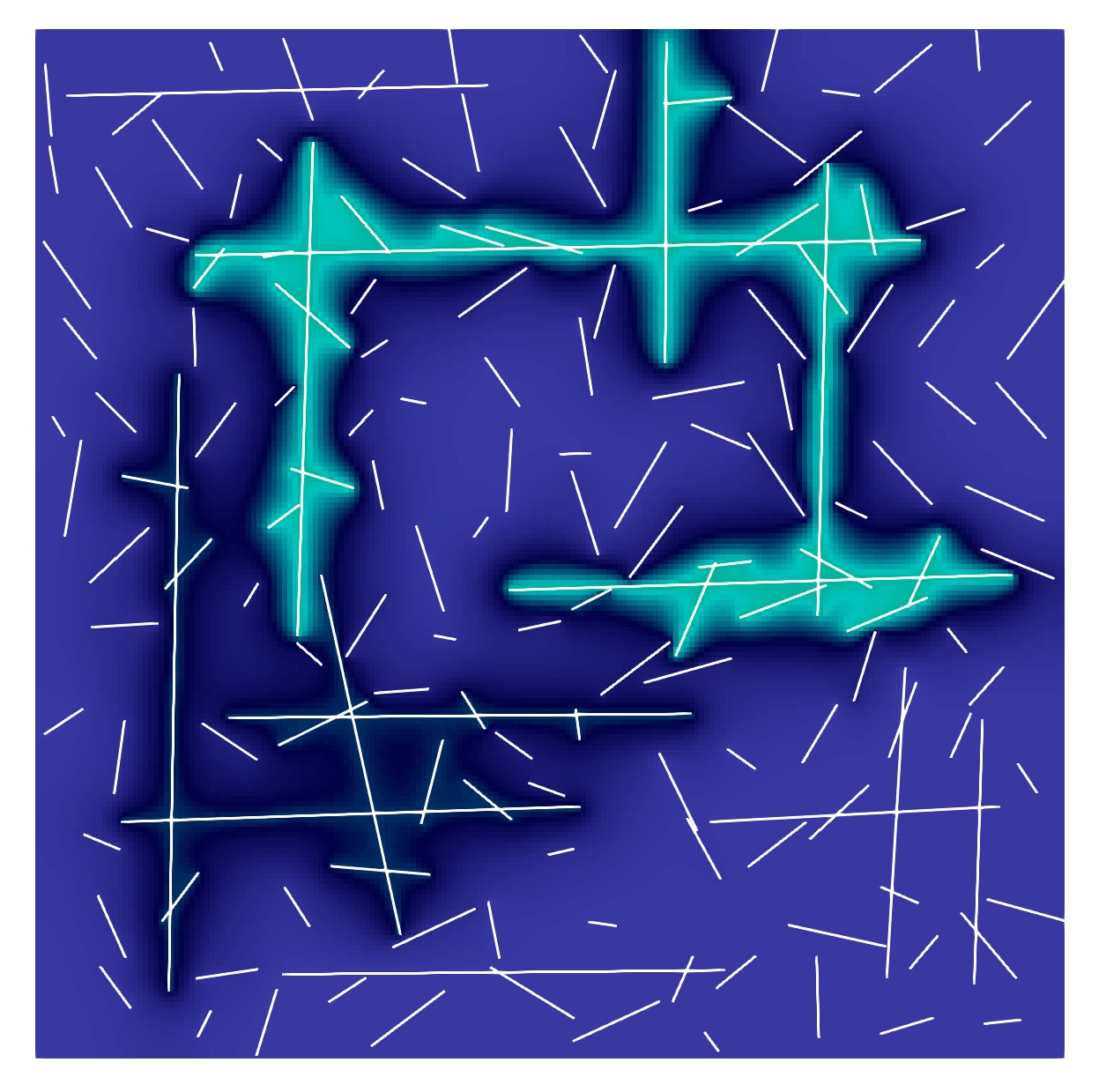}
\includegraphics[width=0.32 \textwidth]{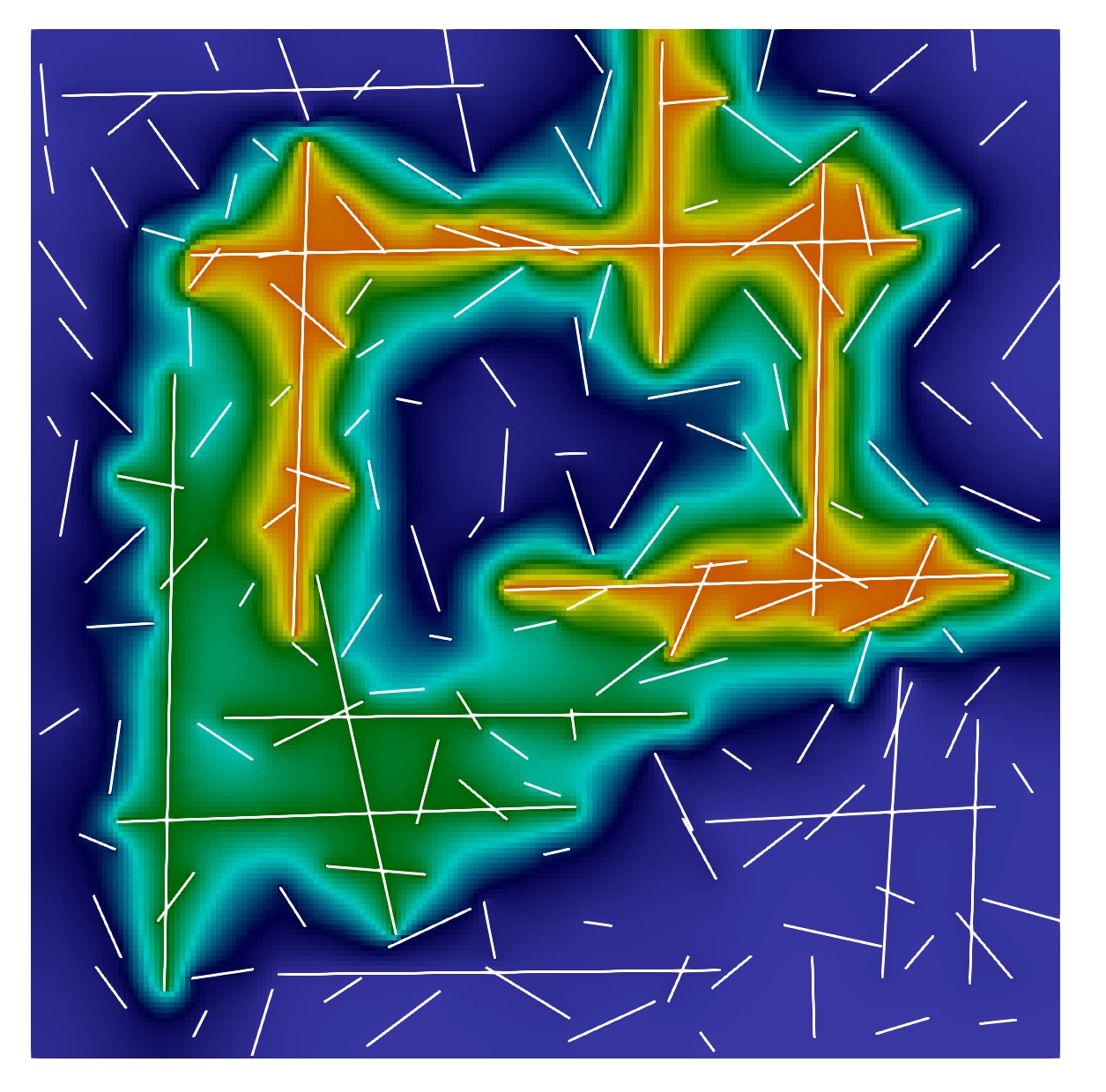}
\includegraphics[width=0.32 \textwidth]{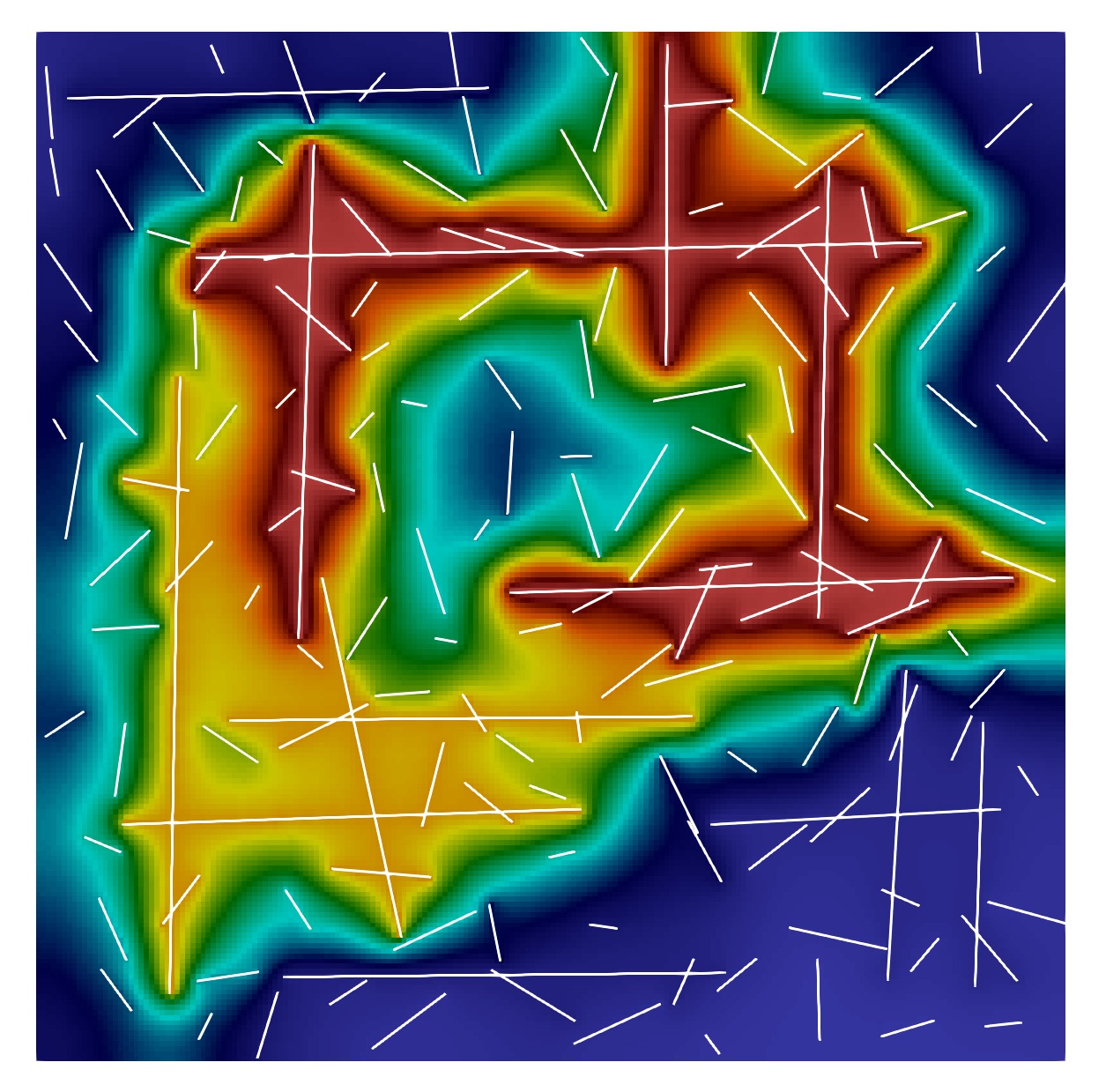}
\caption{Multiscale solutions on mesh $40 \times 40$ with $K^4$ using NLMC model for different time steps $t_{10} = 0.02$,  $t_{30} = 0.06$ and $t_{50} = 0.1$ (from top to bottom). \textit{Geometry 2}. 
First row: upscaled intermediate grid solution. 
Second row: downscaled fine grid solution.  }
\label{fig:u2}
\end{figure}

\begin{table}[h!]
\centering
\begin{tabular}{ |c | c |  c | }
\hline
s &  $e^{FI}_I$ &  $e^{FI}_F$  \\ \hline
1 & 5.466 & 17.626 \\ \hline
2 & 0.416 &  3.917 \\ \hline
3 & 0.112 &  0.901 \\ \hline
4 & 0.103 &  0.236 \\ \hline
6 & 0.101 &  0.104 \\ \hline
\end{tabular}
\,\,\,\,
\begin{tabular}{ |c | c | c | }
\hline
s &  $e^{FI}_I$ &  $e^{FI}_F$  \\ \hline
1 & 50.412 & 51.208  \\ \hline
2 & 1.205  &  4.177 \\ \hline
3 & 0.385 &  0.930 \\ \hline
4 & 0.126 &  0.229 \\ \hline
6 & 0.123 &  0.228 \\ \hline
\end{tabular}
\caption{Relative errors for NLMC intermediate grid solution with different number of oversampling layers $K^s$, $s = 1, 2, 3, 4$ and $6$. 
Left: \textit{Geometry 1} with $DOF_I = 1965$ and $DOF_F = 41042$. 
Right: \textit{Geometry 2} with $DOF_I =2428$ and $DOF_F = 43216$. }
\label{err-nlmc-40}
\end{table}

\begin{figure}[h!]
\centering
\includegraphics[width=0.49 \textwidth]{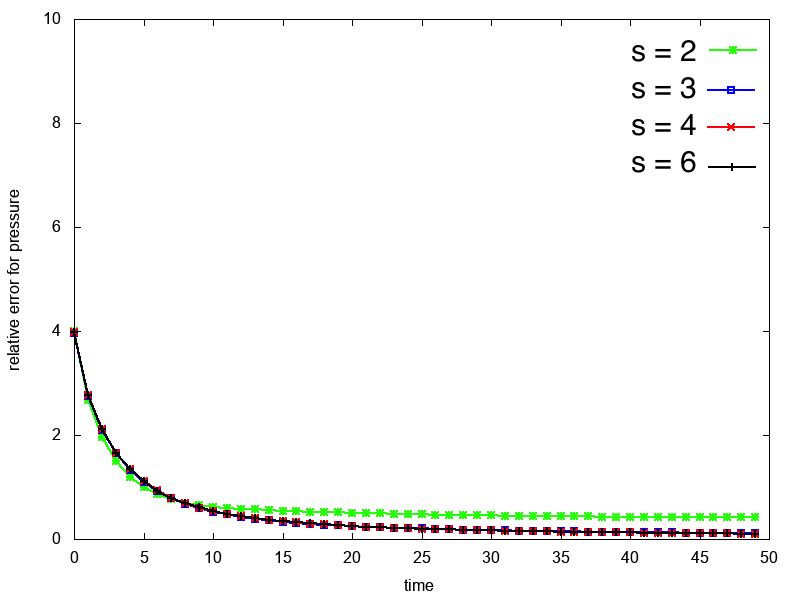}
\includegraphics[width=0.49 \textwidth]{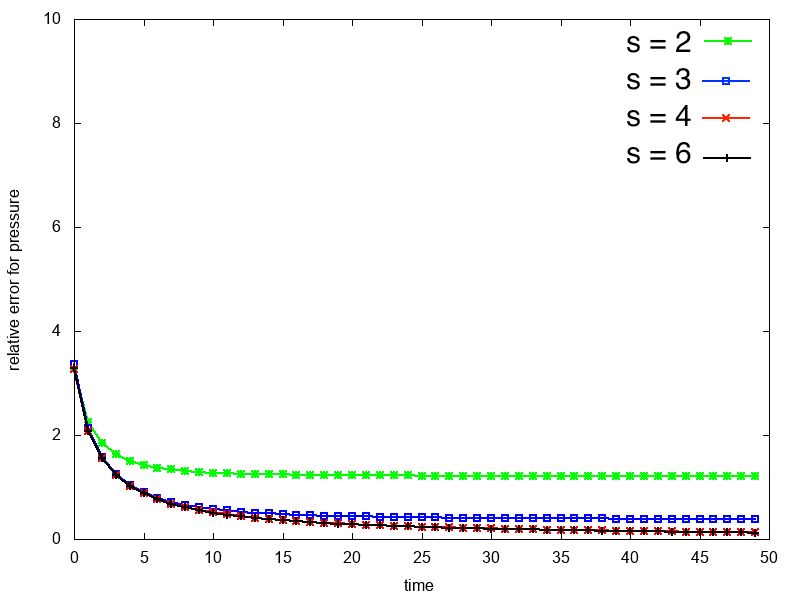}\\
\includegraphics[width=0.49 \textwidth]{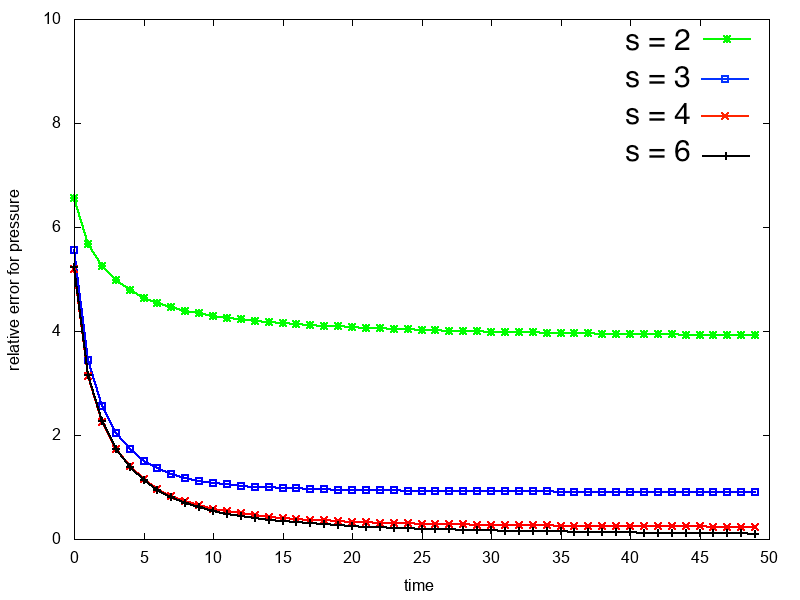}
\includegraphics[width=0.49 \textwidth]{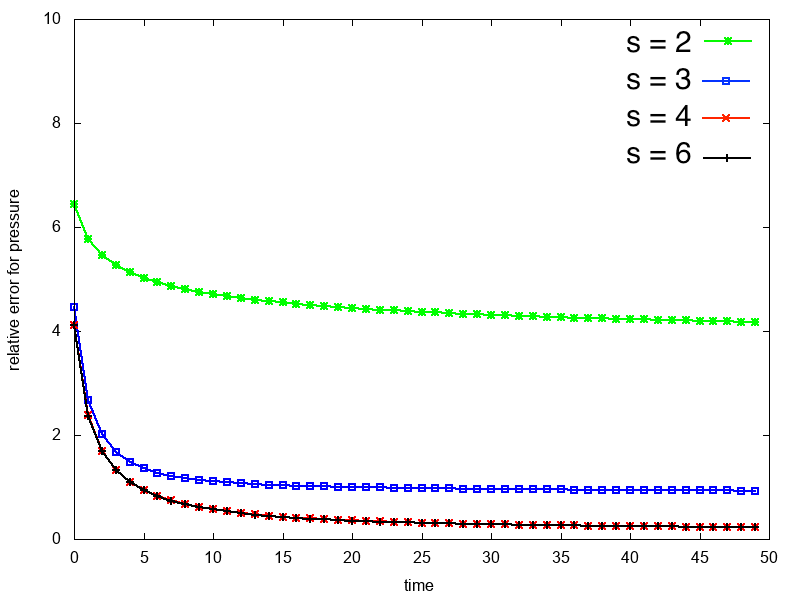}
\caption{Relative errors vs time for upscaled  intermediate grid solution with different number of oversampling layers $K^s$, $s = 2,3,4$ and $6$. 
Left: \textit{Geometry 1}. 
Right: \textit{Geometry 2}. }
\label{tab:t1}
\end{figure}

We construct three grids for multiscale solver:
\begin{itemize}
\item \textit{Fine} level with mesh $200 \times 200$.
\item \textit{Intermediate} level with mesh $40 \times 40$.
\item \textit{Coarse} level with coarse grids $5 \times 5$ and $10 \times 10$.
\end{itemize}
For approximation on fine grid, we constrict finite volume approximation using embedded fracture model. We note that, another approximation techniques can be used, for example, discrete fracture model with unstructured grids.  
Fine grid for fractures domain  for \textit{Geometry 2} contains 3216 cells. For \textit{Geometry 1}, we use grid with 1042 cells for fractures. 
In Figure \ref{fig:mesh}, the fine grid for \textit{Geometry 1} and \textit{Geometry 2} is depicted with blue color and contains $40 000$ cells. The intermediate grid is depicted by red color and contains $1600$ cells. By black color, we depict the coarse grid that contains $36$ and $121$ vertices. 
Note that $DOF_C$, $DOF_I$ and $DOF_F$ are the number of degrees of freedom for coarse, intermediate and fine grids approximations.

We set following parameters for model problem: $a_m = 10^{-5}$, $a_f = 10^{-6}$, $b_m = 10^{-6}$, $b_f = 1.0$ with $\sigma = 10^{-4}$.  We set $p_0 = 0$ as initial pressure and zero flux on boundary.  
We set a source term on the fractures inside cells $K = [0.1, 0.15] \times [0.05, 0.1]$ and $K = [0.6, 0.65] \times [0.9, 0.95]$ with $q = 10^{-3}$.
We simulate $t_{max} = 0.1$ with 50 time steps.

\textbf{Intermediate grid approximation using NLMC method. }

First, we consider relative errors for upscaled multicontinuum model using NLMC method on intermediate grid. 
To compare the results, we use the relative $L^2$ errors between fine grid in upscaled intermediate grid models $e^{FI}$. 
We calculate errors on intermediate grid ($e^{FI}_{I}$) and on fine grid ($e^{FI}_F$)
\[
e^{FI}_{I} = \frac{ || p_I - \bar{p} ||_{L^2} }{ || p_I ||_{L^2} },  \quad
e^{FI}_F = \frac{ || p - \bar{p}_F ||_{L^2}  }{ || p ||_{L^2} },
\]
where $\bar{p}$ is the upscaled intermediate grid solution, $\bar{p}_F = R^T \bar{p}$ is the downscaled of fine grid  intermediate grid solution $\bar{p}$, $p$ is the reference fine grid solution, $p_I$ is the intermediate grid cell average for reference fine grid solution $p$ and 
\[
|| p_I - \bar{p} ||^2_{L^2} =  \sum_K ( p^K_I - \bar{p}^K )^2,
\quad  p^K_I = \frac{1}{|K|} \int_K p \, dx.
\]

In Figures \ref{fig:u1} and \ref{fig:u2}, we present the pressure on mesh $40 \times 40$ with $K^4$ using upscaled model for different time steps $t_{10} = 0.02$,  $t_{30} = 0.06$ and $t_{50} = 0.1$ \textit{Geometry 1} and \textit{Geometry 2}, respectively. 
In the first row, we depict an upscaled medium grid solution. 
Using projection matrix, we can reconstruct fine grid solution from intermediate grid upscaled model (second row in figures).
The fine-scale systems have $DOF_f = 41042$ for \textit{Geometry 1} and $DOF_f = 43216$ for \textit{Geometry 2}. 
Upscaled intermediate grid model has $DOF_c = 1965$ for \textit{Geometry 1} and $DOF_c = 2428$ for \textit{Geometry 2}. 
NLMC method provides accurate meaningful intermediate grid solution with less then one percent errors on fine and intermediate grids.

In Tables \ref{fig:u1} and \ref{fig:u2}, we show relative errors on intermediate and fine grids for different number of oversampling layers $K^s$ with $s = 1,2, 3, 4$ and $6$. 
For intermediate grid approximation with 1600 cells, when we take 4 oversampling layers, we have $0.1 \%$ of intermediate grid error at final time for \textit{Geometry 1} and similar fine grid error. We observe that one oversampling layer cannot provide accurate solution and we should use sufficient number of oversampling layers for obtaining good solution. 
In Figures \ref{tab:t1}, we show relative errors vs time for upscaled  intermediate and fine grids solution with different number of oversampling layers $K^s$, $s = 2, 3, 4$ and $6$.  For intermediate grid solution, we can obtain accurate results with more than 2 oversampling layers. For accurate reconstructed fine grid solution, we should take more than 3 oversampling layers. The proposed method provide accurate solutions for unsteady  mixed dimensional coupled system for fractured porous media for both test geometries and reduce size of the system a lot. For example, we have $DOF_I = 1965$ and $DOF_F = 41042$  for  \textit{Geometry 1}. For  \textit{Geometry 2}, we have $DOF_I =2428$ and $DOF_F = 43216$.

\begin{table}[h!]
\centering
\begin{tabular}{ |c | c | c | c | }
\hline
$M$ &$DOF_C$ &  $e^{IC}_I$ &  $e^{IC}_F$ \\ \hline
1 		& 36 & 49.155  	& 49.475 \\ \hline
4 		& 144 & 9.065  		& 10.146 \\ \hline
8 		& 288 & 7.823  		& 8.917 \\ \hline
12 	& 432 & 4.506  	& 5.210 \\ \hline
16 	& 576 & 2.218  	& 2.634 \\ \hline
20 	& 720 & 1.588  	& 1.903 \\ \hline
24 	& 864 & 0.908  	& 1.116 \\ \hline
28 	& 1008 & 0.370  	& 0.503 \\ \hline
\end{tabular}
\,\,\,\,
\begin{tabular}{ |c | c | c | c | }
\hline
$M$ &$DOF_C$ &  $e^{IC}_I$ &  $e^{IC}_F$ \\ \hline
1 		& 36 & 59.190	& 59.519 \\ \hline
4 		& 144 & 59.189	& 59.518 \\ \hline
8 		& 288 & 58.316	& 58.450 \\ \hline
12 	& 432 & 37.888	& 37.954 \\ \hline
16 	& 576 & 8.046	& 8.417 \\ \hline
20 	& 720 & 3.667	& 4.303 \\ \hline
24 	& 864 & 2.021	& 2.599 \\ \hline
28 	& 1008 & 1.934	& 2.491 \\ \hline
\end{tabular}
\caption{Relative errors for GMsFEM with $5 \times 5$ coarse grid solution with different number of multiscale basis functions $M$. 
Left: \textit{Geometry 1}. Right: \textit{Geometry 2}. }
\label{err-ms-5}
\end{table}

\begin{table}[h!]
\centering
\begin{tabular}{ |c | c | c | c | }
\hline
$M$ &$DOF_C$ &  $e^{IC}_I$ &  $e^{IC}_F$ \\ \hline
1  	 & 121 & 48.473  & 48.616 \\ \hline
2  	 & 242 & 15.437  & 16.173 \\ \hline
4  	 & 484 & 3.949  & 4.531 \\ \hline
8  	 & 968 & 1.177  & 1.446 \\ \hline
12  	 & 1452 & 0.367  & 0.495 \\ \hline
\end{tabular}
\,\,\,\,
\begin{tabular}{ |c | c | c | c | }
\hline
$M$ &$DOF_C$ &  $e^{IC}_I$ &  $e^{IC}_F$ \\ \hline
1  	 & 121 & 59.190  & 59.519 \\ \hline
2  	 & 242 & 59.124  & 59.422 \\ \hline
4  	 & 484 & 42.111  & 41.714 \\ \hline
8  	 & 968 & 3.171  & 3.867 \\ \hline
12  	 & 1452 & 1.336  & 1.772 \\ \hline
\end{tabular}
\caption{Relative errors for GMsFEM with $10 \times 10$ coarse grid.  solution with different number of multiscale basis functions $M$. 
Left: \textit{Geometry 1}. Right: \textit{Geometry 2}. }
\label{err-ms-10}
\end{table}

\begin{figure}[h!]
\centering
\includegraphics[width=0.49 \textwidth]{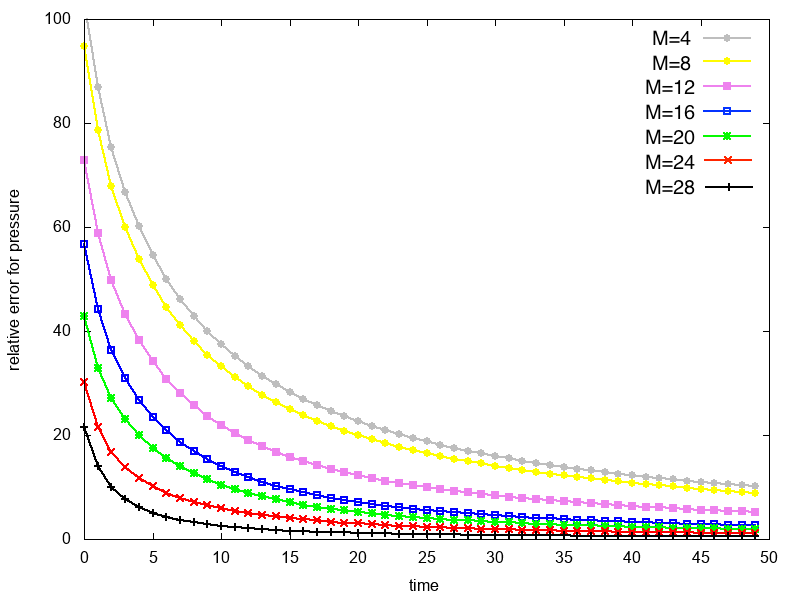}
\includegraphics[width=0.49 \textwidth]{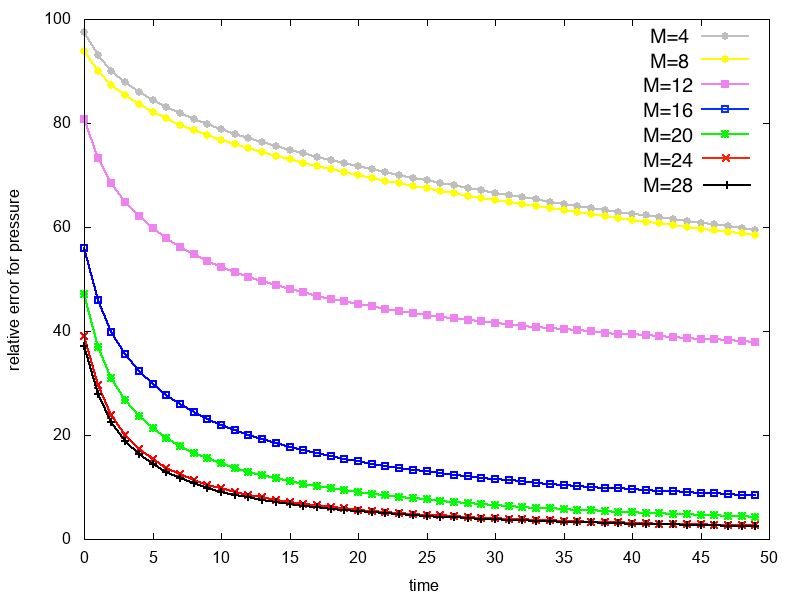}\\
\includegraphics[width=0.49 \textwidth]{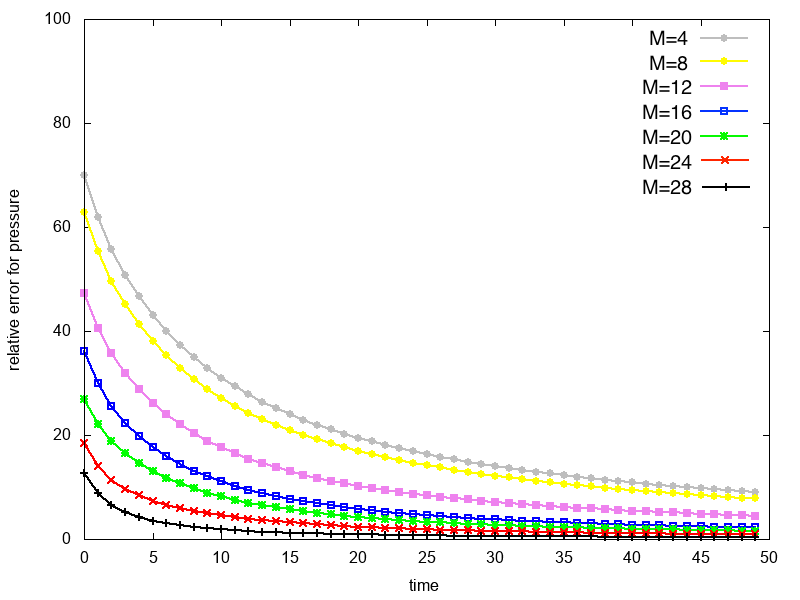}
\includegraphics[width=0.49 \textwidth]{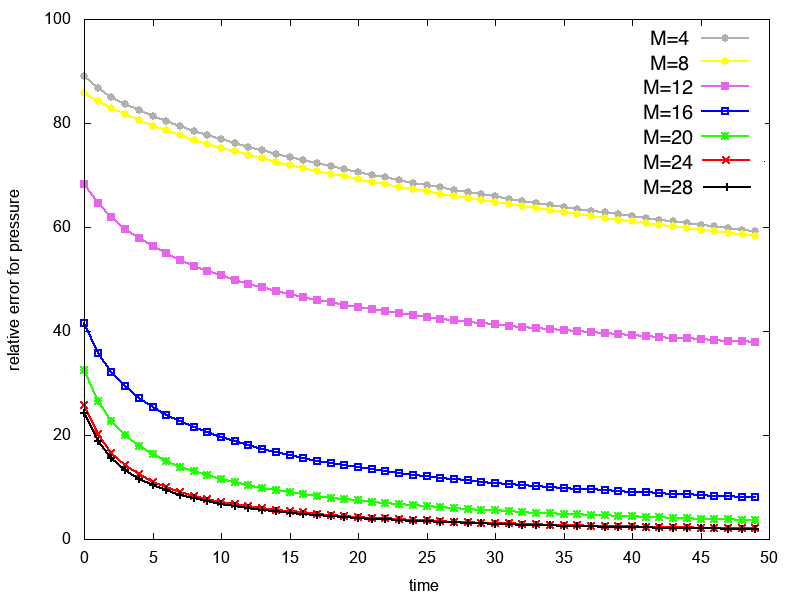}
\caption{Relative errors vs time for GMsFEM with $5 \times 5$ coarse grid.  
First row:  $e^{IC}_I$. Second row:   $e^{IC}_F$.  
Left: \textit{Geometry 1}. 
Right: \textit{Geometry 2}. }
\label{tab:t2}
\end{figure}

\begin{figure}[h!]
\centering
\includegraphics[width=0.49 \textwidth]{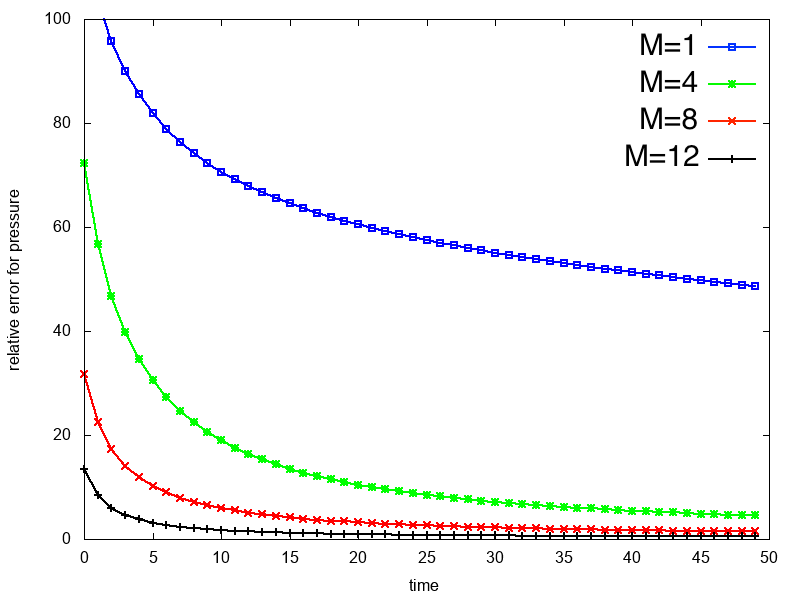}
\includegraphics[width=0.49 \textwidth]{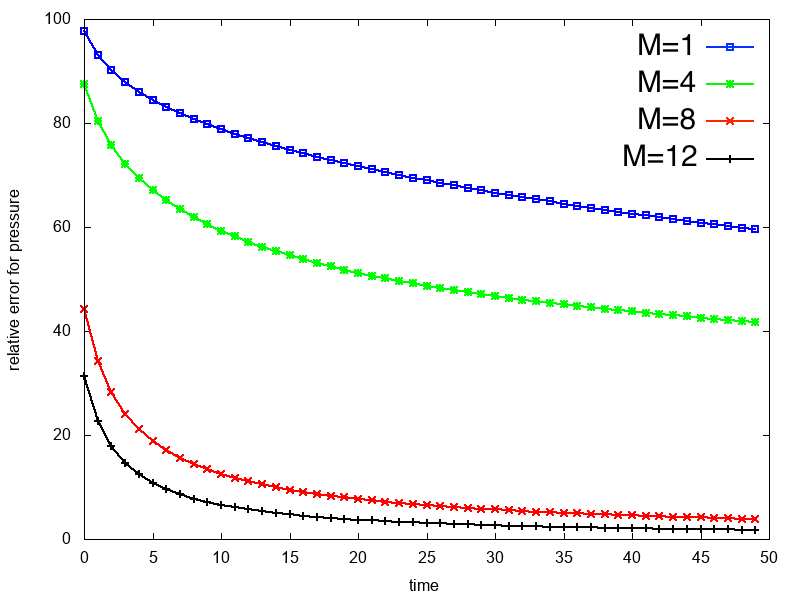}\\
\includegraphics[width=0.49 \textwidth]{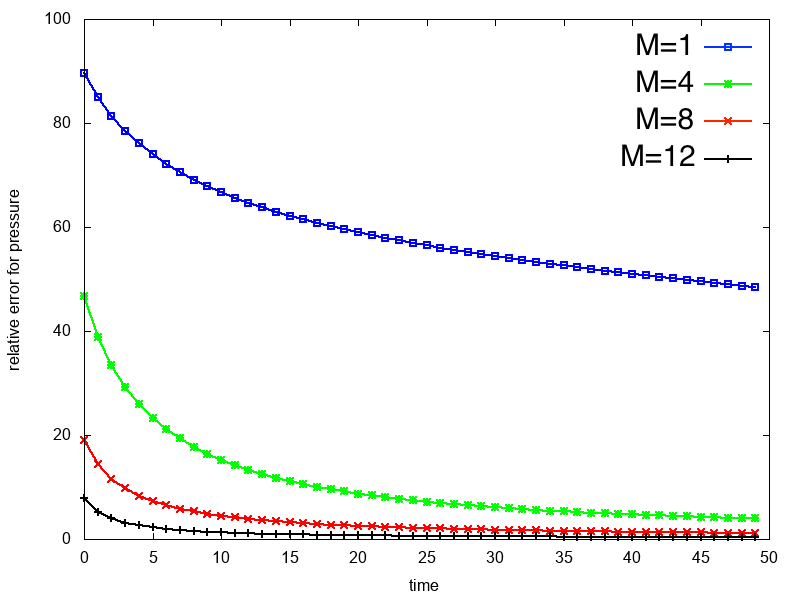}
\includegraphics[width=0.49 \textwidth]{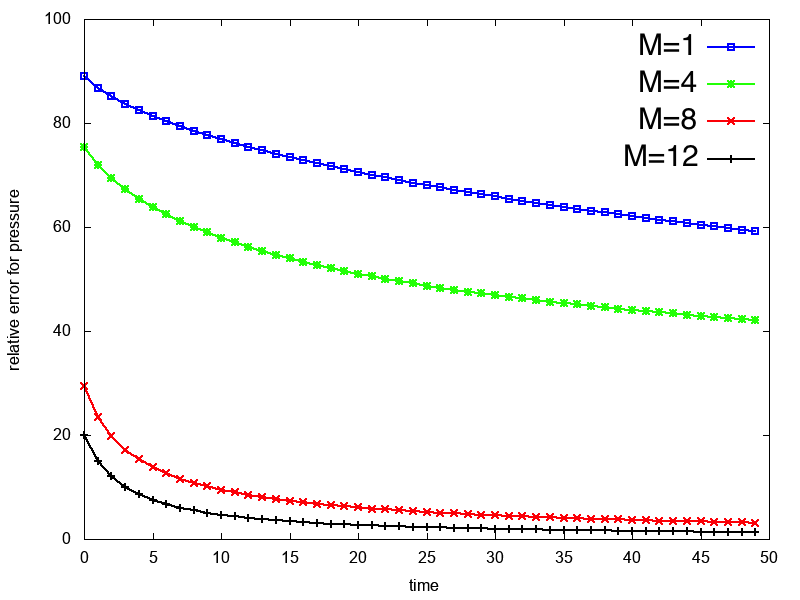}
\caption{Relative errors vs time for GMsFEM with $10 \times 10$ coarse grid. First row:  $e^{IC}_I$. Second row:   $e^{IC}_F$. 
Left: \textit{Geometry 1}. 
Right: \textit{Geometry 2}. }
\label{tab:t3}
\end{figure}

\textbf{Coarse grid approximation using GMsFEM. }

Next, we consider the coarse grid approximation using GMsFEM using intermediate grid upscaled model.  We  use an intermediate grid approximation projection matrix for reconstruction of the fine grid solution. 
We calculate errors between reference fine grid  and GMsFEM solutions on intermediate and fine grids
\[
e^{IC}_I = \frac{ || p_I - p_C ||_{L^2} }{ || p_I ||_{L^2} }, \quad
e^{IC}_F = \frac{ || p - p_F ||_{L^2} }{ || p ||_{L^2} }, 
\]
where $p_C$ is the GMsFEM solution, $p_F = R^T p_C$ is the reconstructed fine grid GMsFEM solution, $p$ is the reference fine grid solution, $p_I$ is the intermediate grid cell average for reference fine grid solution $p$.

We consider two coarse grids: $5 \times 5$ and $10 \times 10$. 
In Tables \ref{err-ms-5} and \ref{err-ms-10}, we show relative errors on intermediate and fine grids for different number of multiscale basis functions, $M$.  The construction of the multiscale basis functions  performed  on intermediate grid. 
For coarse grid approximation with 36 vertices for sufficient number of multiscale basis function, we obtain accurate solution with one percent  of errors for \textit{Geometry 1} and  \textit{Geometry 2}. For finer coarse grid, we can use smaller number of miltiscale basis functions for accurate approximation. 
In Figures \ref{tab:t2} and \ref{tab:t3}, we depict the relative errors vs time for GMsFEM with $5 \times 5$ and $10 \times 10$ coarse grid, respectively. We observe that for geometry with larger number of fractures, we should use more multiscale basis functions.  For example, we obtain $3.9 \%$ of intermediate grid errors  for \textit{Geometry 1}, when we take 4 multiscale basis functions on $10 \times 10$ coarse grids. We obtain similar errors   for \textit{Geometry 2}, when we take 8 multiscale basis functions. 

Finally, we discuss the computational advantages in terms of degrees of freedom. 
In GMsFEM method, we have offline and online stages. On online stage, we calculate multiscale basis functions and construct coarse grid matrices. On offline stage, we solve coarse grid system. 
We can consider proposed algorithm as an extension of the GMsFEM for the upscaled multicontinuum models. 
The  advantage of the proposed method in the acceleration of the GMsFEM model construction by performing offline stage on the intermediate coarse grid for upscaled model, where  nonlocal multicontinuum method used for construction an accurate model. 

Next, we consider the computational advantages of the offline computations. 
Let $N^C_{vert}$ is the number local domains $\omega_i$, $i = 1, ..., N^C_{vert}$. We construct multiscale basis functions in each $\omega$ by solution of the local spectral problems. If we perform calculations of the fine grid, the number of degrees of freedom of local spectral problem is $DOF_{\omega} = N^{\omega}_F$, where for finite volume approximation $N^{\omega}_F = N^{\omega, m}_F + N^{\omega, f}_F$,  $N^{\omega, m}_F$ and $N^{\omega, f}_F$ is the number of fine grid cells for porous matrix and for fractures grid, respectively. When we perform solution on the local spectral problem on intermediate grid using upscaled model, we have $DOF_{\omega} = N^{I, \omega}_{cells}  + \sum_{j = 1}^{N^{I, \omega}_{cells}} L_j $, where $N^{I, \omega}_{cells}$ is the number of intermediate grid cells $K_j$ in $\omega$  and $L^{\omega}_j$ is the number of fractures in $K_j \in \omega$. 
If fine grid is $200 \times 200$ and intermediate grid is $40 \times 40$, then for local domain $\omega_{26}$ and performing calculations on the fine grid,  we have $DOF_{\omega} = 6899$ with $N^{\omega, m}_F = 6400$ and $N^{\omega, f}_F = 499$ for coarse grid $5 \times 5$.
For same coarse grid and same local domain, for the case of intermediate grid based GMsFEM basis construction, we have $DOF_{\omega} = 387$ with $N^{I, \omega}_{cells} = 256$. 
Furthermore, construction of the coarse grid system using intermediate upscaled model can also be done much faster. 
For online computation using GMsFEM on coarse grid $5 \times 5$, we have $DOF_C = 720$ for 20 multiscale basis functions and for fine grid $DOF_F = 41042$ for \textit{Geometry 1}.  

We proposed  three-level technique for multiscale simulations for fractured porous media. On the fine grid we use embedded fracture model, but another methods can be used, for example, discrete fracture model. On intermediate grid, we use nonlocal multicontimuum method to construct an upscaled model. On coarse grid, we construct multiscale solver based on the Generalized Multiscale Finite Element Method. We perform numerical simulations for   three-level  method for  model problems for two fractures geometries.


\bibliographystyle{plain}
\bibliography{lit}

\end{document}